\documentclass[epsfig,latexsym,amsfonts,twoside]{article}
\usepackage{amssymb,enumerate}
\usepackage{amsmath,amscd}

\usepackage{tikz}
\tikzset{help lines/.style={color=blue!50,very thin}}

\def\part#1{\frac{\partial\phantom{#1}}{\partial#1}}
\newtheorem{thm}{Theorem}
\newtheorem{theorem}[thm]{Theorem}
\newtheorem{proposition}[thm]{Proposition}

\newtheorem{conjecture}[thm]{Conjecture}

{\hfill $\Box$ \end{trivlist}}
\newenvironment{definition}{\begin{trivlist}\item[]{\bf Definition}\em }%
{\end{trivlist}}
\newenvironment{remark}{\begin{trivlist}\item[]{\bf Remark} }%
{\end{trivlist}}
\newenvironment{example}{\begin{trivlist}\item[]{\bf Example} }%
{\end{trivlist}}
\newenvironment{question}{\begin{trivlist}\item[]{\bf Question} }%
{\end{trivlist}}

\def\Z{\ifmmode{{\mathbb Z}}\else{${\mathbb Z}$}\fi}
\def\Q{\ifmmode{{\mathbb Q}}\else{${\mathbb Q}$}\fi}
\def\C{\ifmmode{{\mathbb C}}\else{${\mathbb C}$}\fi}
\def\P{\ifmmode{{\mathbb P}}\else{${\mathbb P}$}\fi}

\def\H{\ifmmode{{\mathrm H}}\else{${\mathrm H}$}\fi}

\def\B{\ifmmode{{\mathcal B}}\else{${\mathcal B}$}\fi}
\def\E{\ifmmode{{\mathcal E}}\else{${\mathcal E}$}\fi}
\def\F{\ifmmode{{\mathcal F}}\else{${\mathcal F}$}\fi}
\def\K{\ifmmode{{\mathcal K}}\else{${\mathcal K}$}\fi}
\def\L{\ifmmode{{\mathcal L}}\else{${\mathcal L}$}\fi}
\def\M{\ifmmode{{\mathcal M}}\else{${\mathcal M}$}\fi}
\def\N{\ifmmode{{\mathcal N}}\else{${\mathcal N}$}\fi}
\def\O{\ifmmode{{\mathcal O}}\else{${\mathcal O}$}\fi}
\def\U{\ifmmode{{\mathcal U}}\else{${\mathcal U}$}\fi}
\def\V{\ifmmode{{\mathcal V}}\else{${\mathcal V}$}\fi}
\def\X{\ifmmode{{\mathcal X}}\else{${\mathcal X}$}\fi}

\def\Br{\ifmmode{{\mathrm{Br}}}\else{${\mathrm{Br}}$}\fi}
\def\OG{\ifmmode{\widetilde{\cal M}_4}\else{$\widetilde{\cal M}_4$}\fi}
\def\D{\ifmmode{{\mathcal D}^b}\else{${{\mathcal
    D}^b}$}\fi}
\def\Shah{\ifmmode{\amalg\hspace*{-3.5pt}\amalg}\else{$\amalg\hspace*{-3.5pt}\amalg$}\fi}

\begin{document}

\title{Lagrangian fibrations by Prym varieties}
%\footnote{2010 {\em Mathematics Subject Classification.\/} 14J32, 32Q25.}}
\author{Justin Sawon}
\date{February, 2019}
\maketitle

\begin{abstract}
We survey Lagrangian fibrations of holomorphic symplectic varieties, both compact and non-compact, whose fibres are Jacobians and Prym varieties.
\end{abstract}

\section{Introduction}

In~\cite{hitchin87ii} Hitchin showed that moduli spaces of Higgs bundles are integrable systems. This means that there is a map given by a collection of Poisson-commuting functions from each moduli space of Higgs bundles to a space of half the dimension; here the Poisson structure is the inverse of the natural holomorphic symplectic structure on the moduli space. The general fibres of this map are abelian varieties: Jacobians of curves for $\mathrm{GL}$-Hitchin systems and Prym varieties for other gauge groups.

Moduli spaces of Higgs bundles are non-compact. Turning to compact holomorphic symplectic manifolds, a remarkable result of Matsushita~\cite{matsushita99} (see Theorem~\ref{matsushita}) states that the only fibrations that they can support are ones whose fibres are Lagrangian with respect to the holomorphic symplectic form. Note that by the remark following Theorem~\ref{matsushita}, Lagrangian fibrations and integrable systems are really the same thing in this context. Examples of Lagrangian fibrations on compact holomorphic symplectic manifolds, and orbifolds, include the integrable systems of Beauville-Mukai~\cite{beauville99, mukai84}, Debarre~\cite{debarre99}, Markushevich-Tikhomirov~\cite{mt07}, Arbarello-Sacc{\`a}-Ferretti~\cite{asf15}, and Matteini~\cite{matteini16}, all described in this paper.

 The $\mathrm{GL}$-Hitchin system and the Beauville-Mukai system are both Lagrangian fibrations by Jacobians of curves, so they share similar properties. For instance, they are both isomorphic to their own dual fibrations, at least up to a global twist, because Jacobians are autodual abelian varieites. This duality can be enhanced to a (twisted) Fourier-Mukai transform between the (twisted) derived category of each Lagrangian fibration and the derived category of its dual fibration. Donagi, Ein, and Lazarsfeld~\cite{del97} found a deeper relation between the $\mathrm{GL}$-Hitchin system and the Beauville-Mukai system; namely, the latter can be degenerated to a compactification of the former.
 
 The other Hitchin systems and the other compact Lagrangian fibrations mentioned above are fibrations by Prym varieties. In this article, we explore the relations between these different Lagrangian fibrations. In most case, they are not isomorphic to their dual fibrations. Hausel and Thaddeus~\cite{ht03} showed that the dual of the $\mathrm{SL}$-Hitchin system is the $\mathrm{PGL}$-Hitchin system. Similarly, the dual of the $\mathrm{Sp}(2n,\mathbb{C})$-Hitchin system is the $\mathrm{SO}(2n+1,\mathbb{C})$-Hitchin system. In general, the dual fibration of a Hitchin system is given by taking the Hitchin system whose gauge group is Langlands dual to the original gauge group.
 
 Turning to the compact Lagrangian fibrations by Prym varieties, Menet~\cite{menet14} showed that the dual of a Markushevich-Tikhomirov system is another Markushevich-Tikhomirov system. For the other compact examples, taking their dual fibrations seems to produce new examples of Lagrangian fibrations, making this a potentially fertile direction to explore.
 
 The degeneration of Donagi, Ein, and Lazarsfeld can be generalized to some Lagrangian fibrations by Prym varieties, giving a very concrete connection between certain compact and non-compact examples. In other cases, we can detect some analogies between compact and non-compact examples without there being a clear connection. Still, there remain many non-compact Lagrangian fibrations which do not have obvious compact counterparts (and conversely).

Here is a summary of the paper. In Section~2 we describe the $\mathrm{GL}$-Hitchin system, a non-compact Lagrangian fibration by Jacobians of curves, and how it is related to its dual fibration. In Section~3 we consider compact Lagrangian fibrations, particularly the Beauville-Mukai system. Again we describe how it is related to its dual fibration, we describe the Donagi-Ein-Lazarsfeld degeneration to a compactification of the $\mathrm{GL}$-Hitchin system, and we summarize what is known about the classification of Lagrangian fibrations by Jacobians of curves. In Section~4 we describe the other Hitchin systems, which are fibred by Prym varieties, and the duality relations between them. In Section~5 we turn to compact Lagrangian fibrations by Prym varieties. We describe the Markushevich-Tikhomirov system, its dual fibration, the Matteini system, and some other examples. In Section~6 we mention a few original results of the author (some obtained jointly with Chen Shen): we extend the Donagi-Ein-Lazarsfeld degeneration result to Lagrangian fibrations by Prym varieties, we speculate on the dual fibration of the Matteini system, and we describe the dual fibration of the Debarre system. Finally, in Section~7 we summarize all of the known examples of Lagrangian fibrations by Jacobians and Prym varieties in Table~1; this suggests where we should perhaps search for new examples.

This work was presented at the satellite meeting of the ICM 2018 ``Moduli spaces in algebraic geometry and applications'' held in Campinas, 26--31 July, 2018. The author is grateful to the organizers for the invitation to speak at this meeting, and also grateful for support from the National Science Foundation, award DMS-1555206.

\section{The $\mathrm{GL}$-Hitchin system}

\subsection{Moduli of Higgs bundles}

Higgs bundles on Riemann surfaces and their moduli spaces were introduced by Hitchin~\cite{hitchin87i}. Fix a Riemann surface $\Sigma$ of genus $g\geq 2$.

\begin{definition}
A $\mathrm{GL}(n,\mathbb{C})$-Higgs bundle on $\Sigma$ is a pair $(E,\Phi)$ of a holomorphic bundle $E$ on $\Sigma$ of rank $n$ and a Higgs field
$$\Phi\in \mathrm{H}^0(\Sigma,K\otimes\mathrm{End}E),$$
where $K$ denotes the canonical bundle of $\Sigma$. We say that $(E,\Phi)$ is stable if for all $\Phi$-invariant subbundles $F\subset E$ we have
$$\frac{\mathrm{deg}F}{\mathrm{rank}F}<\frac{\mathrm{deg}E}{\mathrm{rank}E},$$
and semi-stable if instead we have $\leq$.
\end{definition}

Fix the rank $n$ and degree $d$ of $E$, and denote by $\mathcal{M}_{\mathrm{GL}}$ the moduli space of stable Higgs bundles. For coprime $n$ and $d$ this is a smooth quasi-projective variety of dimension $2n^2(g-1)+2$. (If $n$ and $d$ are not coprime, we consider the moduli space of $S$-equivalence classes of semi-stable Higgs bundles; this is a quasi-projective variety, which contains the moduli space of stable Higgs bundles as an open subvariety.)

\begin{thm}[Hitchin~\cite{hitchin87i}]
For coprime $n$ and $d$, $\mathcal{M}_{\mathrm{GL}}$ is a holomorphic symplectic manifold, i.e., it admits a natural holomorphic symplectic form $\sigma$. In fact, $\mathcal{M}_{\mathrm{GL}}$ admits a hyperk{\"a}hler metric.
\end{thm}

\begin{remark}
Denote by ${\mathcal Bun}_{\mathrm{GL}}$ the moduli space of stable bundles on $\Sigma$ of rank $n$ and degree $d$. If $E\in{\mathcal Bun}_{\mathrm{GL}}$, then the inequality above is satisfied for {\em all\/} subbundles $F$. Therefore we could choose any Higgs field $\Phi$, and the pair $(E,\Phi)$ would be a stable Higgs bundle. Consider the cotangent space to ${\mathcal Bun}_{\mathrm{GL}}$ at $E$. By Serre duality
$$T^*_E{\mathcal Bun}_{\mathrm{GL}}=\mathrm{H}^1(\Sigma,\mathrm{End}E)^*\cong \mathrm{H}^0(\Sigma,K\otimes\mathrm{End}E),$$
So in this case $(E,\Phi)$ is really a point of the cotangent bundle $T^*{\mathcal Bun}_{\mathrm{GL}}$. This shows that $T^*{\mathcal Bun}_{\mathrm{GL}}\subset\mathcal{M}_{\mathrm{GL}}$, and in fact it is a dense open subset. The holomorphic symplectic form $\sigma$ on $\mathcal{M}_{\mathrm{GL}}$ is an extension of the canonical holomorphic symplectic form on $T^*{\mathcal Bun}_{\mathrm{GL}}$.
\end{remark}

\subsection{The Hitchin system}

Hitchin~\cite{hitchin87ii} showed that $\mathcal{M}_{\mathrm{GL}}$ admits the structure of an integrable system by defining a map 
\begin{eqnarray*}
h:\mathcal{M}_{\mathrm{GL}} & \longrightarrow & A_{\mathrm{GL}}:=\bigoplus_{i=1}^n\mathrm{H}^0(\Sigma,K^i) \\
(E,\Phi) & \longmapsto & (\mathrm{tr}\Phi,\mathrm{tr}(\Phi^2),\ldots,\mathrm{tr}(\Phi^n))
\end{eqnarray*}
To understand this map, choose a local frame for $E$ and think of the Higgs field $\Phi$ as an $n\times n$ matrix of one-forms. The eigenvalues of this matrix will be $n$ one-forms, and $\mathrm{tr}\Phi,\mathrm{tr}(\Phi^2),\ldots,\mathrm{tr}(\Phi^n)$ are precisely the power sum symmetric polynomials in these one-forms. Although this description is local on $\Sigma$, the unordered set of eigenvalue one-forms do not depend on the choice of frame, and hence their symmetric polynomials are defined globally on $\Sigma$.

Conversely, the symmetric polynomials $\mathrm{tr}\Phi,\mathrm{tr}(\Phi^2),\ldots,\mathrm{tr}(\Phi^n)$ determine the $n$ eigenvalue one-forms locally, giving us an $n$-valued section of $K\rightarrow\Sigma$. In other words, they determine a curve $C$ in the total space of $K$ which maps $n$-to-$1$ to $\Sigma$ under the projection.
$$\begin{array}{ccc}
C & \subset & \mathrm{Tot}K \\
 & _{n:1}\searrow & \downarrow \\
 & & \Sigma \\
\end{array}$$
Moreover, above each eigenvalue we can place its eigenspace, and this produces a line bundle $L$ over $C$, at least for general $(E,\Phi)$. For special $(E,\Phi)$ the curve $C$ could be singular, and $L$ could be a rank-one sheaf that is not locally free.

\begin{definition}
$C$ is called a spectral curve and $(C,L)$ is called the spectral data of $(E,\Phi)$.
\end{definition}

\begin{proposition}
A Higgs bundle $(E,\Phi)$ can be recovered from its spectral data $(C,L)$.
\end{proposition}

Essentially, $E$ is the direct sum of the eigenspaces, which is given by $\pi_*L$ where $\pi:C\rightarrow\Sigma$. To recover $\Phi$ from $(C,L)$, start with the map $L\rightarrow \pi^*K\otimes L$ given by multiplying by the canonical section of $\pi^*K\rightarrow\mathrm{Tot}K$, and take its image under $\pi_*$. The resulting map
$$E=\pi_*L\longrightarrow \pi_*(\pi^*K\otimes L)=K\otimes\pi_*L=K\otimes E$$
is the Higgs field $\Phi$.

\begin{theorem}[Hitchin~\cite{hitchin87ii}]
The map $h:\mathcal{M}_{\mathrm{GL}}\rightarrow A_{\mathrm{GL}}$ makes $\mathcal{M}_{\mathrm{GL}}$ into an algebraic completely integrable system.
\end{theorem}

Note that by Riemann-Roch, $A_{\mathrm{GL}}$ is a vector space of dimension
$$\mathrm{dim}A_{\mathrm{GL}}=\sum_{i=1}^nh^0(\Sigma,K^i)=g+\sum_{i=2}^n(2i-1)(g-1)=n^2(g-1)+1=\frac{1}{2}\mathrm{dim}\mathcal{M}_{\mathrm{GL}}.$$
If we regard $h$ as given by a collection of functions $\{h_i\}$ on $\mathcal{M}_{\mathrm{GL}}$, then the theorem says that these functions Poisson commute,
$$[h_i,h_j]:=\sigma^{-1}(dh_i,dh_j)=0,$$
where the Poisson form $\sigma^{-1}$ is the inverse of the symplectic form $\sigma$ on $\mathcal{M}_{\mathrm{GL}}$.

By general principles of Hamiltonian mechanics, the fibres of an integrable system must be (open subsets of) tori. The Hitchin map $h$ is actually proper, i.e., its fibres are compact. As we saw, a point in the base $A_{\mathrm{GL}}$ determines a spectral curve $C$. The fibre above this point is the Jacobian $\mathrm{Jac}^dC$ for some $d\in\mathbb{Z}$. (Here $d$ is not equal to the degree of the original Higgs bundle, but it can be calculated from it.)

\begin{remark}
There are also singular fibres, corresponding to singular spectral curves. In the mildest cases, the spectral curve $C$ will acquire a node, and the corresponding fibre will be the compactified Jacobian $\overline{\mathrm{Jac}}^dC$, which is the moduli space of rank-one torsion-free sheaves on $C$. In fact, provided $C$ is reduced and irreducible, the corresponding fibre of the Hitchin map will still be the compactified Jacobian $\overline{\mathrm{Jac}}^dC$. But there are also `more singular' fibres, the most singular being the so-called {\em nilpotent cone\/}
$$h^{-1}(0):=\{(E,\Phi)\;|\;\Phi\mbox{ is nilpotent}\}.$$
Here $0\in A_{\mathrm{GL}}$ determines the non-reduced spectral curve $n\Sigma$, i.e., the zero section of $\mathrm{Tot}K$ with multiplicity $n$.
\end{remark}

\subsection{The dual fibration}

Given a fibration by tori, we are often interested in constructing the dual fibration and studying its relation to the original fibration. In mirror symmetry, the Strominger-Yau-Zaslow conjecture states that mirror Calabi-Yau manifolds are dual fibrations in the large complex structure limit. Here the fibrations are by special Lagrangian tori, though note that on hyperk{\"a}hler manifolds there is a deformation of complex structures known as {\em hyperk{\"a}hler rotation\/} that takes holomorphic Lagrangian tori to special Lagrangian tori. 

Mukai~\cite{mukai81} proved that dual abelian varieties are derived equivalent, and this can be extended to the relative setting where we have a family of abelian varieties. Let us apply this to the $\mathrm{GL}$-Hitchin system. Denote by $\mathcal{M}_{\mathrm{GL}}^{sm}$ the union of smooth fibres of the Hitchin system and consider first the case when the fibres are Jacobians $\mathrm{Jac}^0C$ of degree $d=0$. Because $\mathrm{Jac}^0C$ is autodual, the dual fibration $\widehat{\mathcal{M}}_{\mathrm{GL}}^{sm}$ is isomorphic to $\mathcal{M}_{\mathrm{GL}}^{sm}$. In this case we can construct a relative Fourier-Mukai transform between their derived categories
$$\Phi:D^b(\mathcal{M}_{\mathrm{GL}}^{sm})\stackrel{\cong}{\longrightarrow} D^b(\widehat{\mathcal{M}}_{\mathrm{GL}}^{sm})=D^b(\mathcal{M}_{\mathrm{GL}}^{sm})$$
with kernel given by the relative Poincar{\'e} bundle.

\begin{remark}
This will extend over some singular fibres, for example, for the singular fibres given by compactified Jacobians of spectral curves that are reduced and irreducible. We will say more about this in Section~3.3.
\end{remark}

\begin{remark}
For degree $d\neq 0$, we define the dual of $\mathrm{Jac}^dC$ to be the Picard scheme $\mathrm{Pic}^0(\mathrm{Jac}^dC)$. It is isomorphic to $\mathrm{Jac}^0C$. Thus the dual fibration of the degree $d$ $\mathrm{GL}$-Hitchin system $\mathcal{M}_{\mathrm{GL},d}^{sm}$ is isomorphic to the degree $0$ $\mathrm{GL}$-Hitchin system $\mathcal{M}_{\mathrm{GL},0}^{sm}$. In this case there is an obstruction to extending local relative Poincar{\'e} bundles to a global relative Poincar{\'e} bundle, or equivalently, to combining local Fourier-Mukai transforms into a global Fourier-Mukai transform. This obstruction takes the form of a gerbe $\beta$, but the construction above can be modified by twisting by this gerbe and considering derived categories of twisted sheaves, producing
$$\Phi:D^b(\mathcal{M}_{\mathrm{GL},d}^{sm})\stackrel{\cong}{\longrightarrow} D^b(\widehat{\mathcal{M}}_{\mathrm{GL},d}^{sm},\beta)=D^b(\mathcal{M}_{\mathrm{GL},0}^{sm},\beta).$$
These {\em twisted Fourier-Mukai transforms\/} were developed by C{\u a}ld{\u a}raru~\cite{caldararu00}, and applied to Lagrangian fibrations by the author in~\cite{sawon04, sawon08i}.
\end{remark}

\section{Compact systems}

\subsection{Lagrangian fibrations}

In this section we consider compact analogues of the integrable systems of the previous section. In other words, we now assume that the total spaces of the systems are compact.

\begin{definition}
Let $X$ be a compact K{\"a}hler manifold of dimension $2n$. We call $X$ a holomorphic symplectic manifold if it admits a holomorphic two-form $\sigma$ that is non-degenerate in the sense that it induces an isomorphism $T\cong T^*$ (equivalently, the top exterior power $\sigma^{\wedge n}$ trivializes the canonical bundle $K_X=\Omega^{2n}$). Furthermore, we call $X$ an irreducible holomorphic symplectic manifold if it is simply-connected and $\mathrm{H}^0(X,\Omega^2)$ is one-dimensional, generated by $\sigma$.
\end{definition}

The following theorems give very strong restrictions on the structure of fibrations on $X$.

\begin{theorem}[Matsushita~\cite{matsushita99}]
\label{matsushita}
Let $X$ be an irreducible holomorphic symplectic manifold of dimension $2n$ and let $\pi:X\rightarrow B$ be a proper morphism with connected fibres, with $B$ normal and $0<\mathrm{dim}B<2n$. Then
\begin{enumerate}
\item $\mathrm{dim}B=n$,
\item the smooth fibres (and every irreducible component of the singular fibres) of $\pi$ are Lagrangian with respect to $\sigma$,
\item the general fibre is a complex torus.
\end{enumerate}
\end{theorem}

\begin{remark}
By definition, a fibre (or irreducible component of a singular fibre) $F$ is Lagrangian if its tangent space $TF$ is maximal isotropic in $TX$ with respect to the symplectic structure $\sigma$. Dually, $\pi^*T^*B$ is maximal isotropic in $T^*X$ with respect to the Poisson structure $\sigma^{-1}$; but this is exactly the definition of an integrable system. Thus {\em Lagrangian fibrations\/} and integrable systems are really the same thing.

We also see that the fibration is equi-dimensional, as every irreducible component of a fibre must have dimension $n$. 
\end{remark}

\begin{remark}
Note that $X$ is compact and K{\"a}hler but not necessarily projective; nevertheless, a Hodge theoretic argument shows that the fibres of $X\rightarrow B$ are always projective. In particular, the general fibre is an abelian variety.
\end{remark}

\begin{theorem}[Hwang~\cite{hwang08}]
With the same hypotheses as above, if $B$ is smooth then it is isomorphic to $\mathbb{P}^n$.
\end{theorem}

\begin{example}
The Hilbert scheme $S^{[n]}$ of points on a K3 surface $S$ is a crepant resolution of the symmetric product $\mathrm{Sym}^nS$. Beauville proved that $S^{[n]}$ is an irreducible holomorphic symplectic manifold. If $S$ is an elliptic K3 surface, then the map $S\rightarrow\mathbb{P}^1$ induces
$$S^{[n]}\longrightarrow\mathrm{Sym}^nS\longrightarrow\mathrm{Sym}^n\mathbb{P}^1=\mathbb{P}^n,$$
and thus $S^{[n]}$ is a Lagrangian fibration in this case. The (smooth) fibres are products of elliptic curves.
\end{example}

\subsection{The Beauville-Mukai integrable system}

Beauville~\cite{beauville99} and Mukai~\cite{mukai84} discovered and studied another Lagrangian fibration associated to K3 surfaces. Let $S$ be a K3 surface that contains a smooth genus $g$ curve $C$. Riemann-Roch shows that $C$ moves in a $g$-dimensional linear system, $|C|\cong\mathbb{P}^g$. Let $\mathcal{C}/\mathbb{P}^g$ be the universal family of all curves linearly equivalent to $C$. We would like to take the relative compactified Jacobian $\overline{\mathrm{Jac}}^d(\mathcal{C}/\mathbb{P}^g)$ of this family. If the N{\'e}ron-Severi lattice $NS(S)$ of $S$ is generated over $\mathbb{Z}$ by $[C]$ then every curve in the linear system is reduced and irreducible, and therefore its compactified Jacobian is well-defined as the moduli space of rank-one torsion-free sheaves on the curve (see D'Souza~\cite{dsouza79} or Altman and Kleiman~\cite{ak80}). Thus we get a fibration $\overline{\mathrm{Jac}}^d(\mathcal{C}/\mathbb{P}^g)\rightarrow\mathbb{P}^g$ whose general fibre is a $g$-dimensional abelian variety.

More generally, we can choose a polarization $H$ on the K3 surface $S$ and take the moduli space $M(0,[C],1-g+d)$ of $H$-stable sheaves on $S$ with {\em Mukai vector\/} $(0,[C],1-g+d)$; see Mukai~\cite{mukai84}. The general element of this moduli space is again a degree $d$ line bundle on a smooth curve in the linear system $|C|$, thought of as a torsion sheaf on $S$.

\begin{theorem}[Mukai~\cite{mukai84}]
The moduli space $M(0,[C],1-g+d)$ admits a holomorphic symplectic structure. Moreover, if $(0,[C],1-g+d)\in\mathrm{H}^{\bullet}(S,\mathbb{Z})$ is primitive and $H$ is {\em generic\/}, then $M(0,[C],1-g+d)$ is compact and an irreducible holomorphic symplectic manifold. 
\end{theorem}

\begin{remark}
If $NS(S)\cong\mathbb{Z}[C]$ then $(0,[C],1-g+d)$ must be primitive, for any choice of degree $d$, and in this case $M(0,[C],1-g+d)$ coincides with the relative compactified Jacobian $\overline{\mathrm{Jac}}^d(\mathcal{C}/\mathbb{P}^g)$ above. But even when $[C]$ is not primitive, the Mukai vector $(0,[C],1-g+d)$ will still be primitive for some choices of $d$.
\end{remark}

When $X^d:=M(0,[C],1-g+d)$ is an irreducible holomorphic symplectic manifold, the map taking a sheaf $\mathcal{E}$ to its support $\mathrm{Supp}\mathcal{E}\in|C|$ induces a morphism
$$X^d=M(0,[C],1-g+d)\longrightarrow |C|\cong\mathbb{P}^g.$$
By Matsushita's theorem this is a Lagrangian fibration; it is known as the {\em Beauville-Mukai integrable system\/}~\cite{beauville99, mukai84}. In fact, one can see the holomorphic symplectic and Lagrangian structures quite explicitly as follows. Let $t$ be a general point in $\P^g$, let $C_t\subset S$ be the corresponding curve in the linear system $|C|$, and let $X_t$ be the corresponding fibre of $X^d\rightarrow\P^g$. We have a short exact sequence
$$0\longrightarrow TX_t\longrightarrow TX^d|_{X_t}\longrightarrow N_{X_t\subset X^d}\longrightarrow 0.$$
Now $X_t$ is the Jacobian $\mathrm{Jac}^dC_t$ of $C_t$, so $TX_t\cong \H^0(C_t,\Omega^1_{C_t})^*$. On the other hand, $N_{X_t\subset X^d}$ is isomorphic to the pull-back of $T_t\P^g\cong\H^0(C_t,N_{C_t\subset S})$. The latter is isomorphic to $\H^0(C,\Omega^1_{C_t})$, as the short exact sequence
$$0\longrightarrow TC_t\longrightarrow TS|_{C_t}\longrightarrow N_{C_t\subset S}\longrightarrow 0$$
and the triviality of $\mathrm{det} TS=K_S^*$ imply that $N_{C_t\subset S}\cong TC_t^*=\Omega^1_{C_t}$. The holomorphic symplectic form on $TX^d$ comes from the natural pairing between the dual vector spaces $TX_t$ and $N_{X_t\subset X^d}$. It is then clear that $TX_t$ is a Lagrangian subspace.

We can think of the Beauville-Mukai system as a compact analogue of the $\mathrm{GL}$-Hitchin system; we will make the relation between these systems more explicit shortly.

\begin{remark}
If $(0,[C],1-g+d)$ is not primitive, then we can compactify $M(0,[C],1-g+d)$ by adding ($S$-equivalence classes of) semi-stable sheaves. This produces an irreducible holomorphic symplectic variety $M^{ss}(0,[C],1-g+d)$. Only in a few special cases does this singular variety admit a symplectic desingularization (see O'Grady~\cite{ogrady99} and Kaledin, Lehn, and Sorger~\cite{kls06}).
\end{remark}

\subsection{The dual fibration}

As with the $\mathrm{GL}$-Hitchin system, we can try to define the dual fibration of the Beauville-Mukai system. The general fibres are Jacobians and hence autodual, or in the case of arbitrary degree $d$, their dual is
$$\mathrm{Pic}^0(\mathrm{Jac}^dC)\cong\mathrm{Jac}^0C.$$
To what extent can we generalize this isomorphism to singular curves $C$? Remarkably, Esteves, Gagn{\'e}, and Kleiman~\cite{egk02,ek05} proved that if $C$ is reduced and irreducible with at worst nodes and/or cusps as singularities, then
$$\overline{\mathrm{Pic}}^0(\overline{\mathrm{Jac}}^dC)\cong\overline{\mathrm{Jac}}^0C.$$
In particular, this applies to all the curves in the genus $g=2$ Beauville-Mukai system when $[C]$ is primitive, and thus the dual fibration $\widehat{X}^d$ of $X^d$ is isomorphic to $X^0$. This enabled the author to construct a derived equivalence between the derived category of $X^d$ and the twisted derived category of $X^0$. Later, Arinkin extended the above autoduality to compactified Jacobians of arbitrary reduced and irreducible curves $C$ with at worst surficial singularities, i.e., in this case too $\overline{\mathrm{Jac}}^dC$ is irreducible and its dual is
$$\overline{\mathrm{Pic}}^0(\overline{\mathrm{Jac}}^dC)\cong\overline{\mathrm{Jac}}^0C.$$
Moreover, Arinkin extended the derived equivalence to all genus $g$.

\begin{theorem}[Sawon~\cite{sawon08i}, Arinkin~\cite{arinkin13}]
When $[C]$ is primitive the dual fibration of $X^d$ is $X^0$ and there is a twisted Fourier-Mukai transform
$$\Phi:D^b(X^d)\stackrel{\cong}{\longrightarrow} D^b(\widehat{X}^d,\beta)=D^b(X^0,\beta)$$
where $\beta\in\mathrm{H}^2(X^0,\mathcal{O}^*)$ is a gerbe on $X^0$.
\end{theorem}

\subsection{The relation between the Beauville-Mukai and $\mathrm{GL}$-Hitchin systems}

The Beauville-Mukai and $\mathrm{GL}$-Hitchin systems are both relative compactified Jacobians of complete linear systems of curve in symplectic surfaces: in the former case, the curves lie in a K3 surface, and in the latter case, they lie in the total space $\mathrm{Tot}K_{\Sigma}$ of the cotangent bundle of $\Sigma$. By deforming one surface to the other, Donagi, Ein, and Lazarsfeld~\cite{del97} showed that the integrable systems are also related by a deformation. More precisely, one must compactify $\mathrm{Tot}K_{\Sigma}$ and the $\mathrm{GL}$-Hitchin system to make this work. Let us sketch this construction.

Suppose that the curve $\Sigma$ is contained in the K3 surface $S$ and $\mathcal{O}(\Sigma)$ is very ample, giving an embedding
$$S\hookrightarrow\mathbb{P}(\mathrm{H}^0(S,\Sigma)^*)=\mathbb{P}^N.$$
Take the cone over $S$ in $\mathbb{P}^{N+1}$ and intersect it with the pencil of hyperplanes in $\mathbb{P}^{N+1}$ containing $\Sigma$. This gives a pencil of surfaces $\mathcal{S}\rightarrow\mathbb{P}^1$. Moreover
\begin{itemize}
\item the intersection with the hyperplane through the apex of the cone gives $\mathcal{S}_0=\overline{\mathrm{Tot}K_{\Sigma}}$, the one-point compactification of $\mathrm{Tot}K_{\Sigma}$ (this is the cone over the canonical embedding of $\Sigma$ in $\mathbb{P}^{N-1}$),
\item the intersection with any other hyperplane is a surface $\mathcal{S}_t$ isomorphic to $S$.
\end{itemize}
Note that every surface $\mathcal{S}_t$ in this pencil contains a curve isomorphic to $\Sigma$. Curves in the linear system $|n\Sigma|$ in one of the K3 surfaces $\mathcal{S}_t$, $t\neq 0$, will deform to spectral curves in $\mathcal{S}_0$ (or to the additional curves that appear when we compactify $\mathrm{Tot}K_{\Sigma}$). We can now consider the moduli space of stable sheaves on $\mathcal{S}$ that contains degree $d$ line bundles on curves in $|n\Sigma|$ (in either $\mathcal{S}_0$ or any other surface $\mathcal{S}_t$ in the pencil). This moduli space is fibred over $\mathbb{P}^1$.

\begin{theorem}[Donagi-Ein-Lazarsfeld~\cite{del97}]
The fibres of the above fibration over $\mathbb{P}^1$ are Beauville-Mukai systems on $\mathcal{S}_t$ for $t\neq 0$ and the compactification of the $\mathrm{GL}$-Hitchin system on $\Sigma$ induced by the one-point compactification $\mathrm{Tot}K_{\Sigma}\subset \overline{\mathrm{Tot}K_{\Sigma}}$ for $t=0$. In short, the Beauville-Mukai system degenerates to a compactification of the $\mathrm{GL}$-Hitchin system.
\end{theorem}

\begin{remark}
The cone $\overline{\mathrm{Tot}K_{\Sigma}}$ over the canonical embedding of $\Sigma$ is highly singular at its apex. In the above situation, where we assume that $\Sigma$ lies in a K3 surface, $\overline{\mathrm{Tot}K_{\Sigma}}$ can be smoothed to this K3 surface $S$. We call such curves $\Sigma$ {\em K3 curves\/}, and it is an interesting problem to determine when a curve $\Sigma$ is a K3 curve and to reconstruct the K3 surface from the curve; see Arbarello, Bruno, and Sernesi~\cite{abs14, abs17}. Note that the moduli of genus $g$ curves has dimension $3g-3$, whereas K3 curves belong to a family of dimension at most $19+g$ (as algebraic K3 surfaces belong to $19$-dimensional families, and $\mathrm{dim}|\Sigma|=g$). Thus by a parameter count, the general curve of genus $\geq 12$ is {\em not\/} a K3 curve.

If $\Sigma$ is not a K3 curve, then there are obstructions to completely smoothing $\overline{\mathrm{Tot}K_{\Sigma}}$. However, in certain cases one can partially smooth $\overline{\mathrm{Tot}K_{\Sigma}}$ to a surface with a less severe singularity. The resulting integrable systems have not yet been studied.
\end{remark}

\subsection{Classification of Lagrangian fibrations by Jacobians}

It appears that the Beauville-Mukai integrable system is the only (compact) Lagrangian fibration by Jacobians. We make this statement precise in the following conjecture.

\begin{conjecture}
Let $\mathcal{C}\rightarrow\P^n$ be a family of reduced and irreducible curves of arithmetic genus $n$. Suppose that the relative compactified Jacobian $X:=\overline{\mathrm{Jac}}^d(\mathcal{C}/\P^n)$ is an irreducible holomorphic symplectic manifold, and therefore a Lagrangian fibration due to the support map to $\P^n$. Then $X\rightarrow\P^n$ is a Beauville-Mukai system, i.e., the family of curves $\mathcal{C}\rightarrow\P^n$ is a complete linear system of curves on a K3 surfaces.
\end{conjecture}

The conjecture is true in the following cases.

\begin{theorem}[Markushevich~\cite{markushevich96}]
The conjecture is true for genus $n=2$.
\end{theorem}

\begin{theorem}[Sawon~\cite{sawon15ii}]
The conjecture is true for genus $n=3$. If we assume that all the curves in the family $\mathcal{C}\rightarrow\P^n$ are non-hyperelliptic then the conjecture is also true for genus $n=4$ and $5$.
\end{theorem}

\begin{remark}
In the genus $3$ case, one can show that either all of the curves in the family $\mathcal{C}\rightarrow\P^3$ are non-hyperelliptic or all are hyperelliptic. In the latter case, one cannot avoid encountering reducible curves in the family, so the hypothesis of the conjecture must be modified slightly to allow for this.
\end{remark}

\begin{theorem}[Sawon~\cite{sawon15i}]
Denote by $\Delta\subset\P^n$ the discriminant locus parametrizing singular curves in the family $\mathcal{C}\rightarrow\P^n$. Then $\Delta$ is a hypersurface in $\P^n$, and if we assume that the degree of $\Delta$ is at least $4n+2$ then the conjecture is true for all genus $n$.
\end{theorem}
The proof uses {\em coisotropic reduction\/} (see~\cite{sawon09}) to extract a K3 surface from the geometry of the Lagrangian fibration $X\rightarrow\P^n$. This idea was first used by Hurtubise~\cite{hurtubise96} in the local setting, i.e., for a family of curves/Jacobians over a disc.

\section{Other Hitchin systems}

So far we have only considered the Hitchin system on the moduli space of $\mathrm{GL}(n,\mathbb{C})$-Higgs bundles. For other gauge groups the Higgs pair $(E,\Phi)$ must satisfy additional constraints, or admit additional structure. Hitchin~\cite{hitchin87ii} showed that the resulting integrable systems have fibres that are various kinds of Prym varieties. In this section we survey these Hitchin systems for the classical groups. See also Schaposnik's lecture notes~\cite{schaposnik14} for a more detailed account.

\subsection{$\mathrm{SL}$-Hitchin systems}

For $\mathrm{SL}(n,\mathbb{C})$-Higgs bundles $(E,\Phi)$ the determinant bundle $\mathrm{det}E:=\Lambda^nE$ must be trivial and the Higgs field $\Phi$ must be a section of $K\otimes\mathrm{End}_0E$, where $\mathrm{End}_0E$ denotes the trace-free endomorphisms of $E$. As before there is a moduli space $\mathcal{M}_{\mathrm{SL}}$ of stable $\mathrm{SL}$-Higgs bundles and a Hitchin map
\begin{eqnarray*}
h:\mathcal{M}_{\mathrm{SL}} & \longrightarrow & A_{\mathrm{SL}}:=\bigoplus_{i=2}^n\mathrm{H}^0(\Sigma,K^i) \\
(E,\Phi) & \longmapsto & (\mathrm{tr}(\Phi^2),\ldots,\mathrm{tr}(\Phi^n))
\end{eqnarray*}
Note that the $i=1$ term on the right hand side is omitted because $\mathrm{tr}\Phi=0$.

As before there is a spectral curve $C$ in $\mathrm{Tot}K$ determined by $\mathrm{tr}\Phi=0$ and $\mathrm{tr}(\Phi^2),\ldots,\mathrm{tr}(\Phi^n)$, and an $n$-to-$1$ cover $\pi:C\rightarrow\Sigma$. This map induces a norm map on (equivalence classes of) divisors
\begin{eqnarray*}
\mathrm{Nm}:\mathrm{Jac}^dC & \longrightarrow & \mathrm{Jac}^d\Sigma \\
\sum m_jp_j & \longmapsto & \sum m_j\pi(p_j).
\end{eqnarray*}

\begin{definition}
When $d=0$ the Prym variety $\mathrm{Prym}(C/\Sigma)$ of $C\rightarrow\Sigma$ is the connected component of $0$ in $\mathrm{Nm}^{-1}(0)\subset\mathrm{Jac}^0C$. When $d\neq 0$ we can define the Prym variety by fixing a degree $d$ divisor $D$ in $\mathrm{Jac}^d\Sigma$ and taking a connected component of $\mathrm{Nm}^{-1}(D)\subset\mathrm{Jac}^dC$.
\end{definition}

\begin{remark}
All connected components of the fibres of the norm map, for both $d=0$ and $d\neq 0$, are abstractly isomorphic as complex manifolds to $\mathrm{Prym}(C/\Sigma)$, though the isomorphism is not canonical. For $d=0$ we define $\mathrm{Prym}(C/\Sigma)$ to be the connected component of $0$ so that it has a zero and therefore is a complex torus, rather than a principal homogeneous space.
\end{remark}

\begin{remark}
Classically, Prym varieties were defined for unramified double covers $C\rightarrow\Sigma$. In this article we use the term {\em Prym variety\/} far more generally, allowing both $n$-to-$1$ covers for $n>2$ and branched covers. In these more general cases, the Prym variety $\mathrm{Prym}(C/\Sigma)$ is usually not principally polarized.
\end{remark}

Recall that $E$ can be recovered from the spectral data $(C,L)$ as $\pi_*L$. One can show that
$$\Lambda^n\pi_*L\cong\mathrm{Nm}(L)\otimes K^{-n(n-1)/2},$$
and therefore $E$ has trivial determinant bundle if and only if the degree $d$ line bundle $L$ on $C$ lies in $\mathrm{Nm}^{-1}(K^{n(n-1)/2})\subset \mathrm{Jac}^dC$. In other words, the fibres of the $\mathrm{SL}$-Hitchin systems are Prym varieties $\mathrm{Prym}(C/\Sigma)$.

\subsection{$\mathrm{PGL}$-Hitchin systems}

Next consider the gauge group $\mathrm{PGL}(n,\mathbb{C})$. The projective linear group $\mathrm{PGL}(n,\mathbb{C})$ is equal to the quotient $\mathrm{SL}(n,\mathbb{C})/\mathbb{Z}_n$, where $\mathbb{Z}_n$ acts by $n$th roots of unity times the identity matrix. This implies that the moduli space $\mathcal{M}_{\mathrm{PGL}}$ of stable $\mathrm{PGL}$-Higgs bundles is the quotient of $\mathcal{M}_{\mathrm{SL}}$ by the group $\mathrm{Jac}^0\Sigma[n]$ of $n$-torsion points in $\mathrm{Jac}^0\Sigma$, where the line bundle $F\in\mathrm{Jac}^0\Sigma[n]$ acts by $(E,\Phi)\mapsto (E\otimes F,\Phi)$; see Hausel and Thaddeus~\cite{ht03}. Note that the Higgs field $\Phi$ does not change under this action, but we can view it as an endomorphism of $E\otimes F$ because
$$\mathrm{H}^0(\Sigma,K\otimes\mathrm{End}_0(E\otimes F))\cong\mathrm{H}^0(\Sigma,K\otimes\mathrm{End}_0E).$$
In particular, one sees from this description that $\mathcal{M}_{\mathrm{PGL}}$ has orbifold singularities.

As before, we have a Hitchin map
$$h:\mathcal{M}_{\mathrm{PGL}}\longrightarrow A_{\mathrm{PGL}}:=\bigoplus_{i=2}^n\mathrm{H}^0(\Sigma,K^i).$$
To describe the fibres, note that the map $\pi:C\rightarrow\Sigma$ induces the pull-back map $\pi^*:\mathrm{Jac}^0\Sigma\rightarrow\mathrm{Jac}^0C$ (here we just consider the degree $d=0$ case). One can show that the fibre of $h$ over the point corresponding to the spectral curve $C$ is the quotient of $\mathrm{Jac}^0C$ by the action of $\mathrm{Jac}^0\Sigma$ induced by $\pi^*$. Equivalently, the fibre can be described as the quotient of the Prym variety $\mathrm{Prym}(C/\Sigma)$ by the finite group $\mathrm{Jac}^0\Sigma[n]$ (see Hausel and Thaddeus~\cite{ht03}).

Since the $\mathrm{SL}$-Hitchin and $\mathrm{PGL}$-Hitchin systems have the same base $A_{\mathrm{SL}}=A_{\mathrm{PGL}}$, it is instructive to compare their fibres. In fact, we see from the natural defining short exact sequences
\begin{eqnarray*}
\mathrm{SL}: & & 0\longrightarrow \mathrm{Prym}(C/\Sigma)\longrightarrow \mathrm{Jac}^0C\stackrel{\mathrm{Nm}}{\longrightarrow}\mathrm{Jac}^0\Sigma\longrightarrow 0 \\
\mathrm{PGL}: & & 0\longrightarrow\mathrm{Jac}^0\Sigma\longrightarrow\mathrm{Jac}^0C\longrightarrow \widehat{\mathrm{Prym}(C/\Sigma)}\longrightarrow 0,
\end{eqnarray*}
where we used the autoduality of both $\mathrm{Jac}^0C$ and $\mathrm{Jac}^0\Sigma$, that the fibre of the $\mathrm{PGL}$-Hitchin system is actually the dual of the Prym variety $\mathrm{Prym}(C/\Sigma)$. This shows that $\mathcal{M}_{\mathrm{SL}}/A_{\mathrm{SL}}$ and $\mathcal{M}_{\mathrm{PGL}}/A_{\mathrm{PGL}}$ are dual fibrations, at least for smooth fibres.

Based on these observations, and insight from the SYZ conjecture, Hausel and Thaddeus established mirror symmetry for $\mathcal{M}_{\mathrm{SL}}$ and $\mathcal{M}_{\mathrm{PGL}}$. Whereas we just considered the degree $d=0$ case above, they considered non-zero degrees, which leads to gerbes on both sides. 

\begin{theorem}[Hausel-Thaddeus~\cite{ht03}]
The moduli spaces $\mathcal{M}_{\mathrm{SL}}$ and $\mathcal{M}_{\mathrm{PGL}}$ are mirror manifolds. In particular, the stringy Hodge numbers of $\mathcal{M}_{\mathrm{PGL}}$ equal the Hodge numbers of $\mathcal{M}_{\mathrm{SL}}$.
\end{theorem}

\begin{remark}
We need {\em stringy\/} Hodge numbers of $\mathcal{M}_{\mathrm{PGL}}$ because it is an orbifold. On the other hand, $\mathcal{M}_{\mathrm{SL}}$ is smooth. In~\cite{ht03}, Hausel and Thaddeus verified the equality of Hodge numbers for $n=2$ and $3$ only, but the equality was later verified for all $n$ by Groechenig, Wyss, and Ziegler~\cite{gwz17}.
\end{remark}

\subsection{$\mathrm{Sp}$-Hitchin systems}

For $\mathrm{Sp}(2n,\mathbb{C})$-Higgs bundles $(E,\Phi)$ the vector bundle $E$ comes equipped with a non-degenerate skew two-form and the Higgs field $\Phi$ is skew with respect to this form. This means that the eigenvalues of $\Phi$ occur in $\pm$ pairs and the traces of odd powers of $\Phi$ must all vanish. The Hitchin map from the moduli space $\mathcal{M}_{\mathrm{Sp}}$ of stable $\mathrm{Sp}$-Higgs bundles looks like
\begin{eqnarray*}
h:\mathcal{M}_{\mathrm{Sp}} & \longrightarrow & A_{\mathrm{Sp}}:=\bigoplus_{i=1}^n\mathrm{H}^0(\Sigma,K^{2i}) \\
(E,\Phi) & \longmapsto & (\mathrm{tr}(\Phi^2),\mathrm{tr}(\Phi^4),\ldots,\mathrm{tr}(\Phi^{2n}))
\end{eqnarray*}
Because the eigenvalues of $\Phi$ occur in $\pm$ pairs, each spectral curve $C\subset\mathrm{Tot}K$ is invariant under multiplication by $-1$ in the fibres of $K\rightarrow\Sigma$. Note that this action, $\eta\mapsto -\eta$ where $\eta$ is the fibre coordinate, is an involution on $\mathrm{Tot}K$ with fixed locus the zero section. The quotient is the total space $\mathrm{Tot}K^2$ of $K^2\rightarrow\Sigma$. Then $C$ is the branched double cover of a curve $D$ in $\mathrm{Tot}K^2$,
$$\begin{array}{ccc}
C & \subset & \mathrm{Tot}K \\
_{2:1}\downarrow & & _{2:1}\downarrow \\
D & \subset & \mathrm{Tot}K^2.
\end{array}$$

\begin{proposition}
The fibres of the Hitchin map $h:\mathcal{M}_{\mathrm{Sp}}\rightarrow\mathcal{A}_{\mathrm{Sp}}$ are the Prym varieties $\mathrm{Prym}(C/D)$ for the double covers $C\rightarrow D$ described above.
\end{proposition}

\begin{remark}
The polarization type of the abelian variety $\mathrm{Prym}(C/D)$ is $(1,\ldots,1,2,\ldots,2)$ where the number of $2$s is precisely the genus of the curve $D$.
\end{remark}

\subsection{$\mathrm{SO}$-Hitchin systems}

The behavior for the special orthogonal groups is different in the odd and even cases. Here we just consider the gauge group $\mathrm{SO}(2n+1,\mathbb{C})$. For $\mathrm{SO}(2n+1,\mathbb{C})$-Higgs bundles $(E,\Phi)$ the vector bundle $E$ comes equipped with a non-degenerate symmetric two-form and the Higgs field $\Phi$ is skew with respect to this form. The eigenvalues of $\Phi$ consist of $0$ together with $n$ $\pm$ pairs. Once again, the traces of odd powers of $\Phi$ must all vanish and the Hitchin map from the moduli space $\mathcal{M}_{\mathrm{SO}(2n+1,\C)}$ of stable $\mathrm{SO}(2n+1,\C)$-Higgs bundles looks like
\begin{eqnarray*}
h:\mathcal{M}_{\mathrm{SO}(2n+1,\C)} & \longrightarrow & A_{\mathrm{SO}(2n+1,\C)}:=\bigoplus_{i=1}^n\mathrm{H}^0(\Sigma,K^{2i}) \\
(E,\Phi) & \longmapsto & (\mathrm{tr}(\Phi^2),\mathrm{tr}(\Phi^4),\ldots,\mathrm{tr}(\Phi^{2n}))
\end{eqnarray*}
Because $0$ is an eigenvalue of $\Phi$, the zero section of $K\rightarrow\Sigma$ occurs as a component of each spectral curve in $\mathrm{Tot}K$. Discard this component and call the union  of the other components $C$. Then $C\subset\mathrm{Tot}K$ is invariant under multiplication by $-1$ in the fibres of $K\rightarrow\Sigma$, and as before we get a diagram
$$\begin{array}{ccc}
C & \subset & \mathrm{Tot}K \\
_{2:1}\downarrow & & _{2:1}\downarrow \\
D & \subset & \mathrm{Tot}K^2.
\end{array}$$
with $C$ a branched double cover of a curve $D$ in $\mathrm{Tot}K^2$.

\begin{proposition}
The fibres of the Hitchin map $h:\mathcal{M}_{\mathrm{SO(2n+1,\mathbb{C})}}\rightarrow\mathcal{A}_{\mathrm{SO}(2n+1,\mathbb{C})}$ are finite covers of the Prym varieties $\mathrm{Prym}(C/D)$ for the double covers $C\rightarrow D$. In fact, they are the dual abelian varieties $\widehat{\mathrm{Prym}(C/D)}$.
\end{proposition}

\begin{remark}
Earlier we saw that the $\mathrm{SL}$-Hitchin and $\mathrm{PGL}$-Hitchin systems are mirror manifolds, i.e., dual fibrations. Now we see that $\mathcal{M}_{\mathrm{SO}(2n+1,\mathbb{C})}/A_{\mathrm{SO}(2n+1,\mathbb{C})}$ and $\mathcal{M}_{\mathrm{Sp}}/A_{\mathrm{Sp}}$ are also dual fibrations. In general, if we let $\mathrm{^LG}$ denote the Langlands dual group of $G$, then $A_{\mathrm{G}}=A_{\mathrm{^LG}}$ and the dual fibration of $\mathcal{M}_{\mathrm{G}}/A_{\mathrm{G}}$ is $\mathcal{M}_{\mathrm{^LG}}/A_{\mathrm{^LG}}$.
\end{remark}

\section{Compact fibrations by Pryms}

The $\mathrm{Sp}$-Hitchin and $\mathrm{SO}$-Hitchin systems both have spectral curves that are double covers of other curves. These curves are contained in $\mathrm{Tot}K$, respectively, $\mathrm{Tot}K^2$. To construct compact analogues of these integrable systems we should replace the branched double cover $\mathrm{Tot}K\rightarrow\mathrm{Tot}K^2$ by a double cover $S\rightarrow T$ of compact surfaces. In this section we describe various examples that have been studied.

\subsection{The Markushevich-Tikhomirov system}

The first example of this kind was constructed by Markushevich and Tikhomirov~\cite{mt07}. They started with a K3 surface $S$ that is a branched double cover of a degree two del Pezzo surface $T$. In this case $T$ contains an elliptic curve $D$ which is covered by a genus three curve $C$ in $S$. Moreover, $D$ moves in a two-dimensional linear system $|D|\cong\mathbb{P}^2$ and we get a $\mathbb{P}^2$-family of genus three curves $\mathcal{C}$ (an incomplete system of curves linearly equivalent to $C$) in $S$ covering elliptic curves $\mathcal{D}$ in $T$,
$$\begin{array}{ccc}
C & \subset & S \\
_{2:1}\downarrow & & _{2:1}\downarrow \\
D & \subset & T.
\end{array}$$
We can take the Prym varieties of these covers, at least for the smooth curves in this family, and obtain a family of abelian surfaces $\mathrm{Prym}(\mathcal{C}/\mathcal{D})$ over an open subset of $\P^2$.

\begin{theorem}[Markushevich-Tikhomirov~\cite{mt07}]
The relative Prym variety described above admits a natural compactification $\overline{\mathrm{Prym}}(\mathcal{C}/\mathcal{D})$. It is a holomorphic symplectic orbifold of dimension four, and a Lagrangian fibration over $\mathbb{P}^2$.
\end{theorem}

Let us briefly sketch the main ideas behind Markushevich and Tikhomirov's theorem. We defined the Prym variety of $\pi:C\rightarrow D$ as a connected component of a fibre of the norm map $\mathrm{Nm}:\mathrm{Jac}^dC\rightarrow\mathrm{Jac}^dD$. An equivalent definition is as follows. Let $\tau:C\rightarrow C$ be the covering involution and $\tau^*:\mathrm{Jac}^dC\rightarrow\mathrm{Jac}^dC$ the induced involution on the Jacobian. The fixed locus of $\tau^*$ is precisely the abelian subvariety of divisor classes pulled back from $D$. Because $\tau^*$ is an involution, the fixed locus of $-\tau^*$ is a complementary abelian subvariety, and we can define
$$\mathrm{Prym}(C/D):=\mathrm{Fix}^0(-\tau^*),$$
where the superscript $0$ denotes taking a connected component (the connected component containing $0$ in the degree $d=0$ case).

Now $\tau$ also induces a map $\tau^*:|C|\rightarrow |C|$ whose fixed locus is precisely the linear subsystem of curves on $S$ pulled back from curves $D\subset T$. Markushevich and Tikhomirov's approach is to start with the relative compactified Jacobian $\overline{\mathrm{Jac}}^d(\mathcal{C}^{\prime})$ of the family $\mathcal{C}^{\prime}$ of {\em all\/} curves in $S$ linearly equivalent to $C$. This is a Beauville-Mukai integrable system; recall that we define it as a Mukai moduli space of stable sheaves on $S$. Then we would like to define $\overline{\mathrm{Prym}}(\mathcal{C}/\mathcal{D})$ as $\mathrm{Fix}^0(-\tau^*)$ because this picks out $\pi^*|D|$ inside $|C|$ and it picks out the Prym varieties of these curves inside their Jacobians. There is no problem with $\tau^*$: the anti-symplectic involution $\tau$ on $S$ induces an anti-symplectic involution on any Mukai moduli space of stable sheaves on $S$, and in particular
$$\tau^*:\overline{\mathrm{Jac}}^d(\mathcal{C}^{\prime})\longrightarrow\overline{\mathrm{Jac}}^d(\mathcal{C}^{\prime}).$$
What about $-1$? For a smooth curve $C$ it is well defined on $\mathrm{Jac}^0C$ (and even on $\mathrm{Jac}^dC$ if we choose an isomorphism with $\mathrm{Jac}^0C$), where it is given by dualizing a line bundle on $C$, $L\mapsto L^{-1}$. To extend this to all fibres, we can exploit the isomorphism
$$L^{-1}=\mathcal{H}om_C(L,\mathcal{O}_C)\cong\mathcal{E}xt^1_S(L,\mathcal{O}_S(-C))$$
for a line bundle $L$ on a smooth curve $C\subset S$. The last expression is well-defined even for sheaves $L$ in $\overline{\mathrm{Jac}}^0(\mathcal{C}^{\prime})$ supported on singular curves, and so we'd like to define $-1$ to be the involution
\begin{eqnarray*}
\iota:\overline{\mathrm{Jac}}^0(\mathcal{C}^{\prime}) & \longrightarrow & \overline{\mathrm{Jac}}^0(\mathcal{C}^{\prime}) \\
L & \longmapsto & \mathcal{E}xt^1_S(L,\mathcal{O}_S(-C)).
\end{eqnarray*}
The finally difficulty is that $\iota$ preserves $H$-stability only if the polarization $H$ of $S$ is given by $C$ (or a multiple thereof). So we are forced to choose $H=C$. This produces an anti-symplectic involution $\iota$ that commutes with $\tau^*$. Unfortunately the Mukai moduli space $\overline{\mathrm{Jac}}^0(\mathcal{C}^{\prime})$ is singular for this choice of (non-generic) polarization $H=C$.

In any case, we can still proceed. The composition $\iota\circ\tau^*$ is a symplectic involution, and thus its fixed locus is a symplectic subvariety. Markushevich and Tikhomirov define
$$\overline{\mathrm{Prym}}(\mathcal{C}/\mathcal{D}):=\mathrm{Fix}^0(\iota\circ\tau^*)\subset\overline{\mathrm{Jac}}^0(\mathcal{C}^{\prime}),$$
and this is the required holomorphic symplectic orbifold, a Lagrangian fibration over $|D|\cong\mathbb{P}^2$.

\begin{remark}
Markushevich and Tikhomirov actually define three non-isomorphic holomorphic symplectic orbifolds: the first is as above, a second arises from degree $d=2$, and a third is a birational modification of the first given by the Mukai flop of an embedded $\mathbb{P}^2$. They all have $28$ isolated singularities, analytically equivalent to $\mathbb{C}^4/\pm 1$. It is well-known that such singularities do not admit symplectic resolutions.
\end{remark}

\begin{remark}
The fibres of $\overline{\mathrm{Prym}}(\mathcal{C}/\mathcal{D})\rightarrow\mathbb{P}^2$ are abelian surfaces with polarization type $(1,2)$. By comparison, the fibres of the Beauville-Mukai system in dimension four are principally polarized, i.e., of type $(1,1)$, and the fibres of the Debarre system~\cite{debarre99} (described in Section~6.3) in dimension four are of polarization type $(1,3)$.
\end{remark}

\subsection{Dual Prym varieties}

As with fibrations by Jacobians, we wish to describe dual fibrations of fibrations by Prym varieties. But first, in this section, we describe a beautiful geometric construction by Pantazis~\cite{pantazis86} of dual Prym varieties.

This construction applies to branched double covers $C\rightarrow D$ where the curve $D$ is also hyperelliptic. Thus we actually have a tower of branched double covers
$$C\stackrel{2:1}{\longrightarrow}D\stackrel{2:1}{\longrightarrow}\mathbb{P}^1.$$
Pantazis constructs another tower of branched double covers
$$C^{\prime}\stackrel{2:1}{\longrightarrow}D^{\prime}\stackrel{2:1}{\longrightarrow}\mathbb{P}^1$$
as follows. Let $p\in\mathbb{P}^1$ be a general point and suppose that $d_1$ and $d_2$ in $D$ sit above $p$, $c_{11}$ and $c_{12}$ in $C$ sit above $d_1$, and $c_{21}$ and $c_{22}$ in $C$ sit above $d_2$. Above $p$ we define the cover $C^{\prime}$ of $\mathbb{P}^1$ as the pairs of lifts of $\{d_1,d_2\}$. There are four pairs of lifts $\{c_{11},c_{21}\}$, $\{c_{11},c_{22}\}$, $\{c_{12},c_{21}\}$, and $\{c_{12},c_{22}\}$ of $\{d_1,d_2\}$, and therefore $C^{\prime}\rightarrow\mathbb{P}^1$ is a $4$-to-$1$ (branched) cover. In addition, there is an involution on $C^{\prime}$ defined by
$$\{c_{11},c_{21}\}\leftrightarrow \{c_{12},c_{22}\}\qquad\qquad\mbox{and}\qquad\qquad\{c_{11},c_{22}\}\leftrightarrow\{c_{12},c_{21}\}.$$
Quotienting $C^{\prime}$ by this involution gives the curve $D^{\prime}$. Then we have the following geometric description of the dual abelian variety of a Prym variety.

\begin{theorem}[Pantazis's bigonal construction~\cite{pantazis86}]
The abelian varieties $\mathrm{Prym}(C/D)$ and $\mathrm{Prym}(C^{\prime}/D^{\prime})$ are dual.
\end{theorem}

\begin{remark}
Let $D\rightarrow\mathbb{P}^1$ be branched over the points $p_1,\ldots,p_{2s}$ in $\mathbb{P}^1$ and let $C\rightarrow D$ be branched over the points $r_1,\ldots,r_{2t}$ in $D$ whose images in $\mathbb{P}^1$ are $q_1,\ldots,q_{2t}$. In Figure $1$ the $p_i$ are (blue) stars and the $r_j$ and $q_j$ are (red) dots. Then for the second tower of double covers $C^{\prime}\rightarrow D^{\prime}\rightarrow \mathbb{P}^1$ the roles of $p_i$ and $q_j$ are reversed.
\end{remark}

\begin{center}
\begin{tikzpicture}[scale=1.1,>=stealth]
%\draw[step=0.5cm,help lines] (-6,-4) grid (6,4);
\draw (-5.4,2) node {$C$};
\draw (-3,2) ellipse (2.0 and 0.7);
  \draw (-4.3,2.1) .. controls (-4.2,1.9) and (-3.8,1.9) .. (-3.7,2.1);
  \draw (-4.2,2) .. controls (-4.1,2.1) and (-3.9,2.1) .. (-3.8,2);
  \draw (-3.3,2.1) .. controls (-3.2,1.9) and (-2.8,1.9) .. (-2.7,2.1);
  \draw (-3.2,2) .. controls (-3.1,2.1) and (-2.9,2.1) .. (-2.8,2);
  \draw (-2.3,2.1) .. controls (-2.2,1.9) and (-1.8,1.9) .. (-1.7,2.1);
  \draw (-2.2,2) .. controls (-2.1,2.1) and (-1.9,2.1) .. (-1.8,2);

\draw[->] (-3,1.2) -- (-3,0.8);

\draw (-5.4,0) node {$D$};
\draw (-3,0) ellipse (2.0 and 0.7);
  \draw (-3.3,0.1) .. controls (-3.2,-0.1) and (-2.8,-0.1) .. (-2.7,0.1);
  \draw (-3.2,0) .. controls (-3.1,0.1) and (-2.9,0.1) .. (-2.8,0);
  \fill[red] (-4,-0.2) circle[radius=2pt];
  \fill[red] (-3.4,-0.3) circle[radius=2pt];
  \fill[red] (-2.8,-0.3) circle[radius=2pt];
  \fill[red] (-2.2,-0.2) circle[radius=2pt];

\draw[->] (-3,-0.8) -- (-3,-1.2);

\draw (-5.4,-2) node {$\mathbb{P}^1$};
\draw (-3,-2) ellipse (2.0 and 0.7);
  \fill[red] (-4,-2.2) circle[radius=2pt];
  \fill[red] (-3.4,-2.3) circle[radius=2pt];
  \fill[red] (-2.8,-2.3) circle[radius=2pt];
  \fill[red] (-2.2,-2.2) circle[radius=2pt];
  \draw[blue] (-3.8,-1.8) node {$\star$};
  \draw[blue] (-3.2,-1.7) node {$\star$};
  \draw[blue] (-2.6,-1.7) node {$\star$};
  \draw[blue] (-2,-1.8) node {$\star$};

\draw (5.4,2) node {$C^{\prime}$};
\draw (3,2) ellipse (2.0 and 0.7);
  \draw (4.3,2.1) .. controls (4.2,1.9) and (3.8,1.9) .. (3.7,2.1);
  \draw (4.2,2) .. controls (4.1,2.1) and (3.9,2.1) .. (3.8,2);
  \draw (3.3,2.1) .. controls (3.2,1.9) and (2.8,1.9) .. (2.7,2.1);
  \draw (3.2,2) .. controls (3.1,2.1) and (2.9,2.1) .. (2.8,2);
  \draw (2.3,2.1) .. controls (2.2,1.9) and (1.8,1.9) .. (1.7,2.1);
  \draw (2.2,2) .. controls (2.1,2.1) and (1.9,2.1) .. (1.8,2);

\draw[->] (3,1.2) -- (3,0.8);

\draw (5.4,0) node {$D^{\prime}$};
\draw (3,0) ellipse (2.0 and 0.7);
  \draw (3.3,0.1) .. controls (3.2,-0.1) and (2.8,-0.1) .. (2.7,0.1);
  \draw (3.2,0) .. controls (3.1,0.1) and (2.9,0.1) .. (2.8,0);
  \draw[blue] (4,0.2) node {$\star$};
  \draw[blue] (3.4,0.3) node {$\star$};
  \draw[blue] (2.8,0.3) node {$\star$};
  \draw[blue] (2.2,0.2) node {$\star$};

\draw[->] (3,-0.8) -- (3,-1.2);

\draw (5.4,-2) node {$\mathbb{P}^1$};
\draw (3,-2) ellipse (2.0 and 0.7);
  \fill[red] (3.8,-2.2) circle[radius=2pt];
  \fill[red] (3.2,-2.3) circle[radius=2pt];
  \fill[red] (2.6,-2.3) circle[radius=2pt];
  \fill[red] (2,-2.2) circle[radius=2pt];
  \draw[blue] (4,-1.8) node {$\star$};
  \draw[blue] (3.4,-1.7) node {$\star$};
  \draw[blue] (2.8,-1.7) node {$\star$};
  \draw[blue] (2.2,-1.8) node {$\star$};

\end{tikzpicture}\\
\vspace*{3mm}
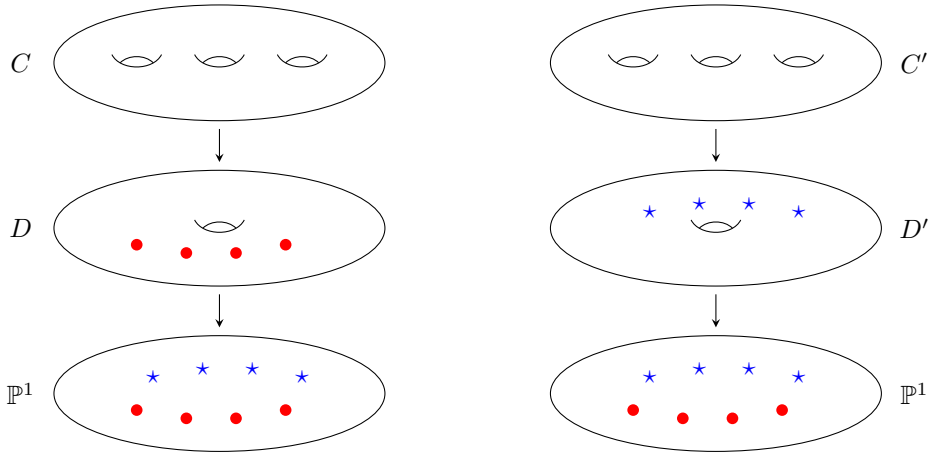
Figure 1: Branch points of towers of double covers
\end{center}

\subsection{Dual fibration of the Markushevich-Tikhomirov system}

Pantazis's bigonal construction can be applied in the relative setting to construct dual fibrations by Prym varieties. To apply this to the Markushevich-Tikhomirov system we first observe that a degree two del Pezzo surface $T$ is a double cover $f:T\rightarrow\mathbb{P}^2$ branched over a plane quartic $\Delta$. Suppose that we have a second plane quartic $\Delta^{\prime}$ that meets $\Delta$ tangentially at exactly eight points. Because of these tangencies, the inverse image $f^{-1}(\Delta^{\prime})$ in $T$ splits into two irreducible components; choose one. Then the double cover of $T$ branched over this component is a K3 surface $S$ and we have a tower of branched double covers
$$S\longrightarrow T\longrightarrow\mathbb{P}^2.$$
(Choosing the other irreducible component of $f^{-1}(\Delta^{\prime})$ gives a K3 surface isomorphic to $S$.) Moreover, every K3 double cover $S$ of the del Pezzo surface $T$ arises in this way, i.e., from a plane quartic $\Delta^{\prime}$ tangent to $\Delta$ at eight points.

As with Pantazis's construction, we can interchange the roles of the two branch loci, $\Delta$ and $\Delta^{\prime}$, to produce a second tower of branched double covers
$$S^{\prime}\longrightarrow T^{\prime}\longrightarrow\mathbb{P}^2.$$
Now the inverse image of a line in $\mathbb{P}^2$ will be an elliptic curve $D$ in $T$ and a genus three curve $C$ in $S$, and similarly an elliptic curve $D^{\prime}$ in $T^{\prime}$ and a genus three curve $C^{\prime}$ in $S^{\prime}$ (as illustrated in Figure $1$). By Pantazis's Theorem the Prym varieties $\mathrm{Prym}(C/D)$ and $\mathrm{Prym}(C^{\prime}/D^{\prime})$ are dual, and relativizing this by varying the line in $\mathbb{P}^2$ produces the following result.

\begin{theorem}[Menet~\cite{menet14}]
The relative compactified Prym varieties $\overline{\mathrm{Prym}}(\mathcal{C}/\mathcal{D})$ and $\overline{\mathrm{Prym}}(\mathcal{C}^{\prime}/\mathcal{D}^{\prime})$ are dual fibrations over $\widehat{\mathbb{P}}^2$. Here $\widehat{\mathbb{P}}^2$ denotes the dual plane which parametrizes lines in $\mathbb{P}^2$.
\end{theorem}
In other words, the dual fibration of a Markushevich-Tikhomirov system is another Markushevich-Tikhomirov system.

\begin{question}
In general the K3 surfaces $S$ and $S^{\prime}$ are not isomorphic, so the dual Markushevich-Tikhomirov systems are not isomorphic. Menet~\cite{menet14} also observes that $S$ and $S^{\prime}$ are not derived equivalent. Are they twisted derived equivalent, i.e., are there gerbes $\beta$ and $\beta^{\prime}$ on $S$ and $S^{\prime}$ such that $D^b(S,\beta)$ is equivalent to $D^b(S^{\prime},\beta^{\prime})$? It is tempting to suspect that the duality between the Markushevich-Tikhomirov systems could be induced by some equivalence like this.
\end{question}

\subsection{The Matteini system and other Prym fibrations}

We describe briefly another Lagrangian fibration by Prym varieties due to Matteini~\cite{matteini16}. It is constructed in a similar way to the Markushevich-Tikhomirov system, except that now we start with a K3 double cover $S\rightarrow T$ of a cubic del Pezzo surface $T$ (degree $3$). In this case we get from
$$\begin{array}{ccc}
C & \subset & S \\
_{2:1}\downarrow & & _{2:1}\downarrow \\
D & \subset & T
\end{array}$$
a $\mathbb{P}^3$-family of genus four curves $\mathcal{C}$ in $S$ covering elliptic curves $\mathcal{D}$ in $T$.
Taking Prym varieties of these covers gives a family of abelian threefolds $\mathrm{Prym}(\mathcal{C}/\mathcal{D})$ over (an open subset of) $\mathbb{P}^3$ which can be compactified in the same way as the Markushevich-Tikhomirov system.

\begin{theorem}[Matteini~\cite{matteini16}]
The relative Prym variety above admits a natural compactification $\overline{\mathrm{Prym}}(\mathcal{C}/\mathcal{D})$ that is a holomorphic symplectic orbifold of dimension six, and a Lagrangian fibration over $\mathbb{P}^3$.
\end{theorem}

\begin{remark}
Matteini describes explicitly the singularities of $\overline{\mathrm{Prym}}(\mathcal{C}/\mathcal{D})$, and notes that they cannot be resolved symplectically.
\end{remark}

\begin{remark}
The fibres of $\overline{\mathrm{Prym}}(\mathcal{C}/\mathcal{D})\rightarrow\mathbb{P}^3$ are abelian threefolds with polarization type $(1,1,2)$.
\end{remark}

In his thesis~\cite{matteini14} Matteini also explored K3 double covers of other del Pezzo surfaces, which yield similar examples to the one above. The common ingredient in these constructions is a K3 surface $S$ with an anti-symplectic involution $\tau$. Nikulin~\cite{nikulin79, nikulin80} gave a compete classification, proving that there are $75$ cases. Then $S$ is a double cover of $T:=S/\tau$, which is necessarily a rational or Enriques surfaces. Some of the examples that have been studied in detail include:
\begin{itemize}
\item the Markushevich-Tikhomirov system~\cite{mt07},
\item the Matteini system~\cite{matteini16},
\item other systems arising from K3 double covers of del Pezzo and Hirzebruch surfaces~\cite{matteini14},
\item the Arbarello-Sacc{\`a}-Ferretti system~\cite{asf15} arising from K3 double covers of Enriques surfaces.
\end{itemize}
In addition to these, one can replace K3 surfaces by abelian surfaces. The Debarre system~\cite{debarre99} (to be described in Section~6.3) is an analogue of the Beauville-Mukai system, constructed from a linear system of curves in an abelian surface $A$. An anti-symplectic involution $\tau$ on $A$ will produce a quotient $T:=A/\tau$ that is a bielliptic surface, and to the above list we can add
\begin{itemize}
\item systems arising from abelian double covers of bielliptic surfaces, also studied by Matteini~\cite{matteini14}.
\end{itemize}

\section{Further directions}

\subsection{Degenerations of Prym fibrations to Hitchin systems}

Earlier we described the Donagi-Ein-Lazarsfeld degeneration of the Beauville-Mukai system to (a compactification of) the $\mathrm{GL}$-Hitchin system. In joint work with Chen Shen, the author has extended this result to some examples whose fibres are Prym varieties.

Recall that the Donagi-Ein-Lazarsfeld degeneration was induced by a degeneration $S\rightsquigarrow \overline{\mathrm{Tot}K_{\Sigma}}$ of a K3 surface to the one-point compactification of $\mathrm{Tot}K_{\Sigma}$. The $\mathrm{Sp}(2n,\mathbb{C})$- and $\mathrm{SO}(2n+1,\mathbb{C})$-Hitchin systems have fibres that are Prym varieties, and finite covers of Prym varieties, respectively. These Prym varieties come from the branched double cover $\mathrm{Tot}K\stackrel{2:1}{\longrightarrow}\mathrm{Tot}K^2$. This suggests taking a K3 surface $S$ that is also a branched double cover, of some surface $T$, and degenerating the double cover as shown.
$$\begin{array}{ccc}
S & \rightsquigarrow & \overline{\mathrm{Tot}K_{\Sigma}} \\
_{2:1}\downarrow & & _{2:1}\downarrow \\
T & \rightsquigarrow & \overline{\mathrm{Tot}K^2_{\Sigma}} \\
\end{array}$$
We can do this in such a way that we get pencils of surfaces $\mathcal{S}\rightarrow\mathbb{P}^1$ and $\mathcal{T}\rightarrow\mathbb{P}^1$ such that
\begin{itemize}
\item $\mathcal{S}_0=\overline{\mathrm{Tot}K_{\Sigma}}$ and $\mathcal{T}_0=\overline{\mathrm{Tot}K^2_{\Sigma}}$,
\item $\mathcal{S}_t\cong S$ and $\mathcal{T}_t\cong T$ for all $t\neq 0$.
\end{itemize}
Then we can construct a relative Prym variety in the usual way, whose fibres are Prym varieties of $C\stackrel{2:1}{\longrightarrow} D$ where $D$ is a curve in $\mathcal{T}_t$ and $C$ its double cover in $\mathcal{S}_t$, for $t\in\mathbb{P}^1$. In other words, we obtain a $\P^1$-family of Prym fibrations, parametrized by $t\in\P^1$.

\begin{theorem}[Chen-Sawon]
The above family gives a degeneration of a Lagrangian fibration by Prym varieties on a compact orbifold (for $t\neq 0$) to a compactification of the $\mathrm{Sp}(2n,\mathbb{C})$-Hitchin system (for $t=0$).
\end{theorem}
Full details will appear elsewhere. Here we just make a few observations.

\begin{remark}
Under the degeneration of surfaces, the branch locus of the double cover $S\rightarrow T$ will degenerate to the branch locus of the double cover $\mathrm{Tot}K_{\Sigma}\rightarrow\mathrm{Tot}K^2_{\Sigma}$. The latter is the zero section of $K^2_{\Sigma}\rightarrow\Sigma$, and hence isomorphic to $\Sigma$ itself. (The double cover $\overline{\mathrm{Tot}K_{\Sigma}}\rightarrow\overline{\mathrm{Tot}K^2_{\Sigma}}$ is also branched over the point at infinity; but this is an additional connected component of the branch locus, and not part of the limit of the branch locus of $S\rightarrow T$.)
\end{remark}

\begin{remark}
The above theorem gives degenerations of {\em some\/} of the Prym fibrations coming from K3 double covers of del Pezzo surfaces to compactifications of $\mathrm{Sp}(2n,\mathbb{C})$-Hitchin systems. However, the Markushevich-Tikhomirov and Matteini systems do not fit into this framework. The reason for this is that in the $\mathrm{Sp}(2n,\mathbb{C})$-Hitchin system the spectral curves lie in the linear system $|n\Sigma|$, where we think of $\Sigma$ as the zero section of $K_{\Sigma}\rightarrow\Sigma$, which is also the branch locus of $\mathrm{Tot}K_{\Sigma}\rightarrow\mathrm{Tot}K^2_{\Sigma}$ (or more precisely, the inverse image of the branch locus). On the other hand, for both the Markushevich-Tikhomirov and Matteini systems the spectral curves $C$ do {\em not\/} lie in the linear system $|nR|$, where $R$ is the branch locus of $S\rightarrow T$ (or more precisely, the inverse image of the branch locus). Since both the spectral curves and the branch loci must be preserved under the degeneration of surfaces, we would get a contradiction.

For example, in the case of the Markushevich-Tikhomirov system the branch locus $R$ comes from half of the pull-back of a plane quartic under $S\rightarrow T\rightarrow\P^2$, whereas $C$ comes from the pull-back of a line. So we would need $n=1/2$, which is impossible.
\end{remark}

\begin{question}
Do the Markushevich-Tikhomirov and Matteini systems degenerate to compactifications of some other Hitchin system?
\end{question}

Recall that the $\mathrm{SO}(2n+1,\C)$-Hitchin system is dual to the $\mathrm{Sp}(2n,\mathbb{C})$-Hitchin system, and it has fibres that are finite covers of Prym varieties.
\begin{question}
Are there Lagrangian fibrations by finite covers of Prym varieties on compact orbifolds (or even manifolds) that degenerate to a compactification of the $\mathrm{SO}(2n+1,\C)$-Hitchin system?
\end{question}

\subsection{Dual fibrations}

As we saw earlier, Menet proved that the dual fibration of a Markushevich-Tikhomirov system is another Markushevich-Tikhomirov system. 

\begin{question}
What is the dual fibration of the Matteini system?
\end{question}
Recall that the fibres of the Matteini system are Prym varieties $\mathrm{Prym}(C/D)$ where $C$ has genus four and $D$ has genus one. These abelian threefolds have polarization of type $(1,1,2)$, and therefore their dual abelian threefolds have polarization of type $(1,2,2)$. We have not yet encountered a Lagrangian fibration with fibres of this polarization type.

In fact, one can say more. Every elliptic curve is a branched double cover of $\P^1$, so we have a tower of branched double covers $C\rightarrow D\rightarrow\P^1$ to which we can apply Pantazis's bigonal construction. The result is another tower of branched double covers $C^{\prime}\rightarrow D^{\prime}\rightarrow\P^1$, where $C^{\prime}$ has genus five and $D^{\prime}$ has genus two, such that the Prym variety $\mathrm{Prym}(C^{\prime}/D^{\prime})$ is dual to $\mathrm{Prym}(C/D)$.

\vspace*{5mm}
\begin{center}
\begin{tikzpicture}[scale=1.1,>=stealth]
%\draw[step=0.5cm,help lines] (-6,-4) grid (6,4);
\draw (-5.4,2) node {$C$};
\draw (-3,2) ellipse (2.0 and 0.7);
  \draw (-4.3,2.1) .. controls (-4.2,1.9) and (-3.8,1.9) .. (-3.7,2.1);
  \draw (-4.2,2) .. controls (-4.1,2.1) and (-3.9,2.1) .. (-3.8,2);
  \draw (-3.3,2.3) .. controls (-3.2,2.1) and (-2.8,2.1) .. (-2.7,2.3);
  \draw (-3.2,2.2) .. controls (-3.1,2.3) and (-2.9,2.3) .. (-2.8,2.2);
  \draw (-2.3,2.1) .. controls (-2.2,1.9) and (-1.8,1.9) .. (-1.7,2.1);
  \draw (-2.2,2) .. controls (-2.1,2.1) and (-1.9,2.1) .. (-1.8,2);
  \draw (-3.3,1.9) .. controls (-3.2,1.7) and (-2.8,1.7) .. (-2.7,1.9);
  \draw (-3.2,1.8) .. controls (-3.1,1.9) and (-2.9,1.9) .. (-2.8,1.8);

\draw[->] (-3,1.2) -- (-3,0.8);

\draw (-5.4,0) node {$D$};
\draw (-3,0) ellipse (2.0 and 0.7);
  \draw (-3.3,-0.1) .. controls (-3.2,-0.3) and (-2.8,-0.3) .. (-2.7,-0.1);
  \draw (-3.2,-0.2) .. controls (-3.1,-0.1) and (-2.9,-0.1) .. (-2.8,-0.2);
  \draw[blue] (-4.2,0.2) node {$\star$};
  \draw[blue] (-3.7,0.3) node {$\star$};
  \draw[blue] (-3.2,0.35) node {$\star$};
  \draw[blue] (-2.7,0.35) node {$\star$};
  \draw[blue] (-2.2,0.3) node {$\star$};
  \draw[blue] (-1.7,0.2) node {$\star$};

\draw[->] (-3,-0.8) -- (-3,-1.2);

\draw (-5.4,-2) node {$\mathbb{P}^1$};
\draw (-3,-2) ellipse (2.0 and 0.7);
  \fill[red] (-4,-2.2) circle[radius=2pt];
  \fill[red] (-3.4,-2.3) circle[radius=2pt];
  \fill[red] (-2.8,-2.3) circle[radius=2pt];
  \fill[red] (-2.2,-2.2) circle[radius=2pt];
  \draw[blue] (-4.2,-1.8) node {$\star$};
  \draw[blue] (-3.7,-1.7) node {$\star$};
  \draw[blue] (-3.2,-1.65) node {$\star$};
  \draw[blue] (-2.7,-1.65) node {$\star$};
  \draw[blue] (-2.2,-1.7) node {$\star$};
  \draw[blue] (-1.7,-1.8) node {$\star$};

\draw (5.4,2) node {$C^{\prime}$};
\draw (3,2) ellipse (2.0 and 0.7);
  \draw (4.3,2) .. controls (4.2,1.8) and (3.8,1.8) .. (3.7,2);
  \draw (4.2,1.9) .. controls (4.1,2) and (3.9,2) .. (3.8,1.9);
  \draw (3.3,1.9) .. controls (3.2,1.7) and (2.8,1.7) .. (2.7,1.9);
  \draw (3.2,1.8) .. controls (3.1,1.9) and (2.9,1.9) .. (2.8,1.8);
  \draw (2.3,2) .. controls (2.2,1.8) and (1.8,1.8) .. (1.7,2);
  \draw (2.2,1.9) .. controls (2.1,2) and (1.9,2) .. (1.8,1.9);
  \draw (3.8,2.3) .. controls (3.7,2.1) and (3.3,2.1) .. (3.2,2.3);
  \draw (3.7,2.2) .. controls (3.6,2.3) and (3.4,2.3) .. (3.3,2.2);
  \draw (2.8,2.3) .. controls (2.7,2.1) and (2.3,2.1) .. (2.2,2.3);
  \draw (2.7,2.2) .. controls (2.6,2.3) and (2.4,2.3) .. (2.3,2.2);

\draw[->] (3,1.2) -- (3,0.8);

\draw (5.4,0) node {$D^{\prime}$};
\draw (3,0) ellipse (2.0 and 0.7);
  \draw (3.8,0.3) .. controls (3.7,0.1) and (3.3,0.1) .. (3.2,0.3);
  \draw (3.7,0.2) .. controls (3.6,0.3) and (3.4,0.3) .. (3.3,0.2);
  \draw (2.8,0.3) .. controls (2.7,0.1) and (2.3,0.1) .. (2.2,0.3);
  \draw (2.7,0.2) .. controls (2.6,0.3) and (2.4,0.3) .. (2.3,0.2);
  \fill[red] (3.8,-0.2) circle[radius=2pt];
  \fill[red] (3.2,-0.3) circle[radius=2pt];
  \fill[red] (2.6,-0.3) circle[radius=2pt];
  \fill[red] (2,-0.2) circle[radius=2pt];

\draw[->] (3,-0.8) -- (3,-1.2);

\draw (5.4,-2) node {$\mathbb{P}^1$};
\draw (3,-2) ellipse (2.0 and 0.7);
  \fill[red] (3.8,-2.2) circle[radius=2pt];
  \fill[red] (3.2,-2.3) circle[radius=2pt];
  \fill[red] (2.6,-2.3) circle[radius=2pt];
  \fill[red] (2,-2.2) circle[radius=2pt];
  \draw[blue] (4.3,-1.8) node {$\star$};
  \draw[blue] (3.8,-1.7) node {$\star$};
  \draw[blue] (3.3,-1.65) node {$\star$};
  \draw[blue] (2.8,-1.65) node {$\star$};
  \draw[blue] (2.3,-1.7) node {$\star$};
  \draw[blue] (1.8,-1.8) node {$\star$};

\end{tikzpicture}\\
\vspace*{3mm}
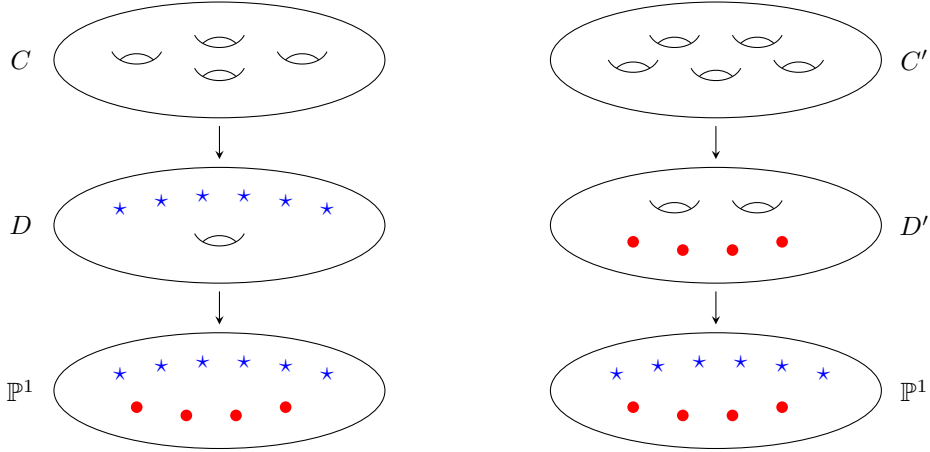
Figure 2: Towers of double covers yielding dual Prym varieties
\end{center}

The next step is to try to combine all of these genus five curves $C^{\prime}$ and genus two curves $D^{\prime}$ into linear systems on a K3 surface $S^{\prime}$ and on a del Pezzo surface $T^{\prime}$, respectively, so that we can build a Lagrangian fibration by Prym varieties. (Note that from the $\P^3$-families of curves $\mathcal{C}$ and $\mathcal{D}$ we get, by Pantazis's construction, $\P^3$-families of curves $\mathcal{C}^{\prime}$ and $\mathcal{D}^{\prime}$; but a priori, there is no guarantee that $\mathcal{C}^{\prime}/\P^3$ and $\mathcal{D}^{\prime}/\P^3$ are linear systems of curves on surfaces.)

In joint work with Chen Shen, the author has approached this from a different direction. We start with a K3 surface $S^{\prime}$ that is a branched double cover of a degree one del Pezzo surface $T^{\prime}$. Then $T^{\prime}$ contains a genus two curve $D^{\prime}$ covered by a genus five curve $C^{\prime}$ in $S^{\prime}$. Indeed, $D^{\prime}$ can be taken to be a general element of the linear system $|-2K_T|$, which is three-dimensional. Thus we have a $\P^3$-family of genus five curves $\mathcal{C}^{\prime}$ in $S^{\prime}$ covering genus two curves $\mathcal{D}^{\prime}$ in $T^{\prime}$.

\begin{proposition}
The relative compactified Prym variety $\overline{\mathrm{Prym}}(\mathcal{C}^{\prime}/\mathcal{D}^{\prime})$ can be constructed as before, i.e., as the fixed locus of a symplectic involution on the relative compactified Jacobian $\overline{\mathrm{Jac}}^0(\mathcal{C}^{\dagger})$, where $\mathcal{C}^{\dagger}$ is the family of {\em all\/} curves in $S^{\prime}$ linearly equivalent to $C^{\prime}$. It is a holomorphic symplectic orbifold of dimension six, and a Lagrangian fibration over $\P^3$.
\end{proposition}
The previous discussion suggests that this is the dual fibration of a Matteini system. However, given a K3 double cover $S^{\prime}$ of a degree one del Pezzo surface $T^{\prime}$, it is not immediately clear how to construct a `dual' K3 double cover $S$ of a degree three del Pezzo surface $T$, or conversely. In fact, the following parameter counts reveals that the situation cannot be this simple:
\begin{itemize}
\item there is a four-dimensional family of degree three del Pezzo surfaces $T$, and the branch locus $\Delta$ of the cover $S\rightarrow T$ belongs to a nine-dimensional linear system $|-2K_T|$; this yields a $13$-dimensional family of K3 double covers of degree three del Pezzo surfaces, $S\rightarrow T$, and thus a $13$-dimensional family of Matteini systems,
\item there is an eight-dimensional family of degree one del Pezzo surfaces $T^{\prime}$, and the branch locus $\Delta^{\prime}$ of the cover $S^{\prime}\rightarrow T^{\prime}$ belongs to a three-dimensional linear system $|-2K_{T^{\prime}}|$; this yields an $11$-dimensional family of K3 double covers of degree one del Pezzo surfaces, $S^{\prime}\rightarrow T^{\prime}$, and thus an $11$-dimensional family of Lagrangian fibrations $\overline{\mathrm{Prym}}(\mathcal{C}^{\prime}/\mathcal{D}^{\prime})$ as in the proposition above.
\end{itemize}

\begin{conjecture}
The dual fibration of the Lagrangian fibration $\overline{\mathrm{Prym}}(\mathcal{C}^{\prime}/\mathcal{D}^{\prime})$ is a special Matteini system, or possibly a degeneration of Matteini systems.
\end{conjecture}
This still leaves open the following question.
\begin{question}
What is the dual fibration of a general Matteini system? Does it admit a description as a relative Prym variety of a linear system of curves on some surface?
\end{question}

\subsection{The Debarre system and its dual fibration}

We have described compact analogues of $\mathrm{Sp}(2n,\mathbb{C})$-Hitchin systems. These examples have fibres that are Prym varieties of covers $C\rightarrow D$, i.e., kernels of surjective maps $\mathrm{Jac}^0C\rightarrow\mathrm{Jac}^0D$, where both $C$ and $D$ vary in a family. By contrast, the $\mathrm{SL}$-Hitchin system has fibres that are kernels of surjective maps $\mathrm{Jac}^0C\rightarrow\mathrm{Jac}^0\Sigma$, where the spectral curves $C$ vary in a family but the curve $\Sigma$, and its Jacobian $\mathrm{Jac}^0\Sigma$, is fixed. In this section we describe a compact Lagrangian fibration by Prym varieties that exhibits analogous behaviour; the construction of this fibration is due to Debarre~\cite{debarre99}.

Just as the dual fibration of the $\mathrm{SL}$-Hitchin system is the $\mathrm{PGL}$-Hitchin system, the dual fibration of the Debarre system can be described in a similar way. We give this description and then end with a conjectural relationship, inspired by mirror symmetry, between the (stringy) Hodge numbers of the Debarre system and its dual fibration.

To construct the Debarre system we begin with an abelian surface $A$ that contains a smooth curve $C$ of genus $g+2$. Moreover, let us assume that the N{\'e}ron-Severi lattice $NS(A)$ of $A$ is generated over $\mathbb{Z}$ by the class $[C]$. Then $C$ is ample and it is a polarization of type $(1,g+1)$ on $A$, as $[C]$ is primitive. Riemann-Roch shows that $C$ moves in a $g$-dimensional linear system, $|C|\cong\P^g$. Let $\mathcal{C}/\mathbb{P}^g$ be the family of curves linearly equivalent to $C$. The assumption on $NS(A)$ ensures that every curve in the linear system is reduced and irreducible, and therefore the relative compactified Jacobian
$$Y:=\overline{\mathrm{Jac}}^d(\mathcal{C}/\mathbb{P}^g)\longrightarrow\mathbb{P}^g$$
is well-defined (see Altman and Kleiman~\cite{ak80}). We could instead define $Y$ as a Mukai moduli space~\cite{mukai84} of $H$-stable sheaves on the abelian surface $A$ with Mukai vector $(0,[C],1-g+d)$, where $H=C$. This yields the same space, but it reveals the holomorphic symplectic structure. However, $Y$ is not an {\em irreducible\/} holomorphic symplectic manifold: it is {\em reducible\/} in the sense that, after taking a finite {\'e}tale cover, it splits into a product of lower dimensional manifolds. Another way of looking at this is that the Albanese map
$$\mathrm{Alb}:Y\longrightarrow \mathrm{Alb}(Y)\cong A$$
is a locally trivial fibration. Up to a finite {\'e}tale cover, $Y$ is the product of $A$ and a fibre of $\mathrm{Alb}$. It is the latter that we are really interested in.

\begin{theorem}[Debarre~\cite{debarre99}]
The fibre of the Albanese map,
$$X:=\mathrm{Alb}^{-1}(0)\subset Y,$$
is an irreducible holomorphic symplectic manifold. Moreover, the restriction of the projection $Y\rightarrow\P^g$ to $X\rightarrow\P^g$ makes $X$ into a Lagrangian fibration.
\end{theorem}

\begin{remark}
Debarre showed that $X$ is deformation equivalent to the generalized Kummer variety $K_g(A)$ of dimension $2g$, defined by Beauville~\cite{beauville83} as the fibre of the natural map
$$\mathrm{Hilb}^{g+1}A\longrightarrow\mathrm{Sym}^{g+1}A\longrightarrow A.$$
In fact, if $d=0$ our Lagrangian fibration is also {\em birational\/} to the generalized Kummer variety $K_g(\widehat{A})$ constructed from the dual abelian surface $\widehat{A}$, as can be seen by applying the methods of Yoshioka~\cite{yoshioka01}.
\end{remark}

\begin{remark}
Let $Y_t=\mathrm{Jac}^dC$ be a fibre of $Y\rightarrow\P^g$ (for some curve $C$ in the linear system on $A$). The restriction of the Albanese map is a surjective homomorphism $Y_t\rightarrow A$, and thus the fibre $X_t$ of $X\rightarrow\P^g$ fits into an exact sequence
$$0\longrightarrow X_t\longrightarrow Y_t=\mathrm{Jac}^dC\longrightarrow A\longrightarrow 0.$$
Since $\mathrm{Jac}^dC$ is principally polarized, $X_t$ must have polarization of type $(1,\ldots,1,g+1)$, complementary to the polarization type $(1,g+1)$ of $A$.
\end{remark}

Next we consider the dual fibration. Dualizing the above exact sequence gives
$$0\longrightarrow \widehat{A}\longrightarrow \widehat{Y}_t=\widehat{\mathrm{Jac}^dC}\longrightarrow \widehat{X}_t\longrightarrow 0.$$
Now the Jacobian is autodual; more canonically, $\widehat{Y}_t=\widehat{\mathrm{Jac}^dC}=\mathrm{Jac}^0C$. Let us restrict to the case $d=0$ so that $\widehat{Y}_t=Y_t$. As we have seen, this autoduality holds also for compactified Jacobians of singular curves, provided the curve are reduced, irreducible, and with surficial singularities (all of which apply here). The conclusion is that $\widehat{X}_t$ is the quotient of $Y_t$ by the action of $\widehat{A}$ induced by the inclusion $\widehat{A}\rightarrow\widehat{Y}_t=Y_t$. This action has a very natural description: $M\in\widehat{A}=\mathrm{Pic}^0A$ is a line bundle on $A$, and it acts on $L\in Y_t=\mathrm{Jac}^0C$ by
$$L\longmapsto M|_C\otimes L.$$
This action clearly extends to singular curves $C$ too, where $L$ may be a rank-one torsion-free sheaf, not necessarily locally free.

In summary, we have a natural action of $\widehat{A}$ on $Y$, and the dual fibration $\widehat{X}$ is the quotient $Y/\widehat{A}$. A priori this would be an algebraic stack, but we can show that it is a smooth Deligne-Mumford stack.

\begin{theorem}[Sawon]
The dual fibration $\widehat{X}:=Y/\widehat{A}$ is a holomorphic symplectic orbifold.
\end{theorem}
For the proof, one shows that the action of $\widehat{A}$ on $Y$ has finite stabilizers. The details will appear elsewhere.

\begin{remark}
The fibres of the dual fibration $\widehat{X}\rightarrow\P^g$ have polarization of type $(1,g+1,\ldots,g+1)$, which is dual to the polarization type $(1,\ldots,1,g+1)$ of the fibres of the Debarre system $X\rightarrow\P^g$. So for $g\geq 3$ we obtain a new example of a Lagrangian fibration on a holomorphic symplectic orbifold.

When $g=2$ both the Debarre system $X\rightarrow\P^2$ and its dual $\widehat{X}\rightarrow\P^2$ have fibres of polarization type $(1,3)$. If an abelian variety is not principally polarized then it will only be isogenous to its dual, not isomorphic, so $X\rightarrow\P^2$ is certainly not self-dual. Moreover, it appears that $\widehat{X}$ is a genuine orbifold, i.e., not smooth, so it is not isomorphic to a Debarre system for any choice of abelian surface $A$.
\end{remark}

Finally we come to the mirror symmetry relation between the Debarre system and its dual fibration. Recall that the fibres of the $\mathrm{SL}$- and $\mathrm{PGL}$-Hitchin systems are dual Prym varieties that appear in dual short exact sequences:
\begin{eqnarray*}
\mathcal{M}_{\mathrm{SL}}: & & 0\longrightarrow \mathrm{Prym}(C/\Sigma)\longrightarrow \mathrm{Jac}^0C\longrightarrow\mathrm{Jac}^0\Sigma\longrightarrow 0 \\
\mathcal{M}_{\mathrm{PGL}}: & & 0\longrightarrow\mathrm{Jac}^0\Sigma\longrightarrow\mathrm{Jac}^0C\longrightarrow \widehat{\mathrm{Prym}(C/\Sigma)}\longrightarrow 0
\end{eqnarray*}
Here we just consider the degree $d=0$ case, and the Jacobians are autodual. Hausel and Thaddeus~\cite{ht03} showed that the stringy Hodge numbers of $\mathcal{M}_{\mathrm{PGL}}$ are equal to the Hodge numbers of $\mathcal{M}_{\mathrm{SL}}$.

For the fibres of the Debarre system and its dual fibration we have (again for degree $d=0$):
\begin{eqnarray*}
X: & & 0\longrightarrow X_t\longrightarrow \mathrm{Jac}^0C\longrightarrow A\longrightarrow 0 \\
\widehat{X}: & & 0\longrightarrow\widehat{A}\longrightarrow\mathrm{Jac}^0C\longrightarrow \widehat{X}_t\longrightarrow 0
\end{eqnarray*}
The analogy with the $\mathrm{SL}$- and $\mathrm{PGL}$-Hitchin systems suggests the following.

\begin{conjecture}
The stringy Hodge numbers of $\widehat{X}$ equal the Hodge numbers of $X$.
\end{conjecture}
As we mentioned earlier, the Debarre system $X$ is deformation equivalent to the generalized Kummer variety $K_g(A)$, whose Hodge numbers have been calculated by G{\"o}ttsche and Soergel~\cite{gs93}. To prove the conjecture one therefore needs to calculate the stringy Hodge numbers of the dual fibration $\widehat{X}$. We can do this in the (rather trivial) $g=1$ case, as follows.

\begin{example}
When $g=1$, $X$ is a (smooth) K3 surface, and an elliptic fibration over $\P^1$. The general fibre of $X\rightarrow\P^1$ is the kernel of a surjective map $\mathrm{Jac}^0C\rightarrow A$ where $C$ is a smooth genus three curve; this of course gives a smooth elliptic curve. The singular fibres of $X\rightarrow\P^1$ correspond to nodal curves $C$. A calculation shows that there are precisely twelve nodal curves in the pencil $\mathcal{C}\rightarrow\P^1$ of genus three curves on $A$, and therefore $X\rightarrow\P^1$ has twelve singular fibres. These fibres sit inside the compactified Jacobians $\overline{\mathrm{Jac}}^0C$ and look like the red curve in the right hand picture of Figure 3; in other words, they are singular elliptic curves of Kodaira type $I_2$.

\begin{center}
\vspace*{3mm}
\begin{tikzpicture}[scale=1.1,>=stealth]
%  \draw[step=0.5cm,help lines] (-7,-2) grid (6,3);
%  \draw (-2,0) -- (-6,0) -- (-6,2.5) -- (-2,2.5) -- (-2,0);
  \draw (-2,0) -- (-2,2.5);
  \draw (-6,0) -- (-6,2.5);
  \draw[thick, blue] (-2,0) -- (-6,0);
  \draw[thick, blue] (-2,2.5) -- (-6,2.5);
  \draw[very thick, red] (-5,0) -- (-5,2.5);
  \draw[very thick, red] (-3,0) -- (-3,2.5);
  \draw[very thin, gray, ->] (-4.9,2.4) -- (-3.1,0.1);
  \draw[very thin, gray] (-2.9,2.4) -- (-2.1,1.3);
  \draw[very thin, gray, ->] (-5.9,1.2) -- (-5.1,0.1);
  
  \draw[->] (-1.2,1.2) -- (-0.7,1.2);
  
  \draw (0,0) .. controls (0.2,0.1) and (0.7,0.5) .. (0.8,0.6);
  \draw[dashed] (0.8,0.6) .. controls (1.5,1.3) and (1.5,2.6) .. (0.8,2.7);
  \draw (0.8,2.7) .. controls (0.1,2.5) and (0.1,1.2) .. (0.8,0.6);
  \draw[dashed] (0.8,0.6) .. controls (0.7,0.7) and (1.2,0.3) .. (1.4,0.2);
  
  \draw (4,-1) .. controls (4.2,-0.9) and (4.7,-0.5) .. (4.8,-0.4);
  \draw (4.8,-0.4) .. controls (5.5,0.3) and (5.5,1.6) .. (4.8,1.7);
  \draw (4.8,1.7) .. controls (4.1,1.5) and (4.1,0.2) .. (4.8,-0.4);
  \draw (4.8,-0.4) .. controls (4.7,-0.3) and (5.2,-0.7) .. (5.4,-0.8);
  
  \draw (0,0) -- (4,-1);
  \draw (0.8,2.7) -- (4.8,1.7);
  \draw[thick, blue] (0.8,0.6) -- (4.8,-0.4);
  \draw[dashed] (1.4,0.2) -- (4.6,-0.6);
  \draw (4.6,-0.6) -- (5.4,-0.8);

  \draw[very thick, red] (0.5,-0.125) .. controls (1,0) and (1.2,0) .. (1.8,0.35);
  \draw[dashed, very thick, red] (1.8,0.35) .. controls (2.4,0.7) and (3.3,2.0) .. (2.8,2.2);
  \draw[very thick, red] (2.8,2.2) .. controls (2.0,1.9) and (3.2,0.15) .. (3.8,-0.15);
  \draw[dashed, very thick, red] (3.8,-0.15) .. controls (4.3,-0.45) and (4.3,-0.5) .. (4.7,-0.6);
  
  \draw[very thick, red] (2.5,-0.625) .. controls (3,-0.5) and (3.2,-0.5) .. (3.8,-0.15);
  \draw[dashed, very thick, red] (3.8,-0.15) .. controls (4.4,0.2) and (5.3,1.5) .. (4.8,1.7);
  \draw[very thick, red] (0.8,2.7) .. controls (0.0,2.4) and (1.2,0.65) .. (1.8,0.35);
  \draw[dashed, very thick, red] (1.8,0.35) .. controls (2.3,0.05) and (2.3,0) .. (2.7,-0.1);
  
\end{tikzpicture}\\
\vspace*{5mm}
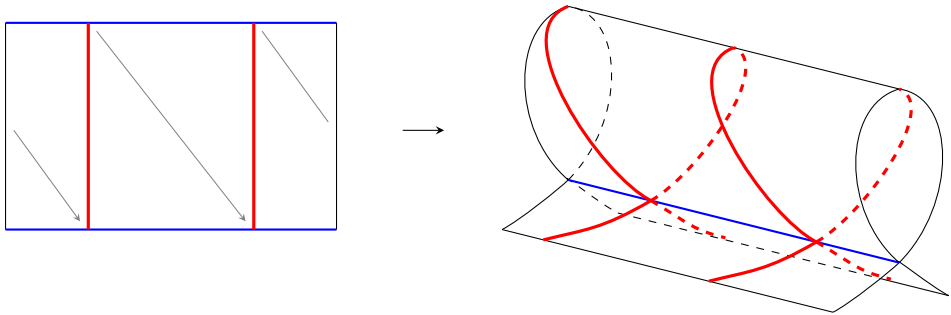
Figure 3: Compactified Jacobian leading to a singular Prym variety of type $I_2$ 
\end{center}

The dual fibration $\widehat{X}\rightarrow\P^1$ must necessarily have twelve singular fibres too. One can show that these look like singular elliptic curves of Kodaira type $I_1$, i.e., nodal rational curves, but the surface $\widehat{X}$ is now a singular K3 surface with twelve $A_1$ singularities at the nodes of these singular fibres. Blowing up these singularities gives an elliptic K3 surface $\widetilde{X}$ with twelve $I_2$ fibres. (Indeed $\widetilde{X}\cong X$, but we do not need this fact.) In this situation, the stringy Hodge numbers of $\widehat{X}$ will equal the Hodge numbers of its crepant resolution $\widetilde{X}$, which are the same as the Hodge numbers of $X$ because $\widetilde{X}$ and $X$ are both smooth K3 surfaces. This proves the conjecture for $g=1$.
\end{example}

\section{Summary of Lagrangian fibrations}

In Table 1 we summarize the different Lagrangian fibrations that we have described in this article.

\begin{table}[h!]
  \begin{center}
    \begin{tabular}{|c|c|c|}  \hline     
    \textbf{Fibres} & \textbf{Non-compact} & \textbf{Compact} \\
      \hline\hline
      Jacobians & $\mathrm{GL}(n,\mathbb{C})$-Hitchin system, \S2.2 & Beauville-Mukai system, \S3.2 \\
      \hline
      \mbox{Prym varieties} & $\mathrm{SL}(n,\mathbb{C})$-Hitchin system, \S4.1 & Debarre system, \S6.3 \\
       & $\mathrm{PGL}(n,\mathbb{C})$-Hitchin system, \S4.2 & dual Debarre system, \S6.3 \\
       & $\mathrm{Sp}(2n,\mathbb{C})$-Hitchin system, \S4.3 & \S6.1 \\
       & $\mathrm{SO}(2n+1,\mathbb{C})$-Hitchin system, \S4.4 & ? \\
       & $\mathrm{SO}(2n,\mathbb{C})$-Hitchin system, \S4.4 & ? \\
       & ? & Markushevich-Tikhomirov system, \S5.1 \\
       & ? & Matteini system, \S5.4 \\
       & ? & dual Matteini system, \S6.2 \\
       & ? & Arbarello-Sacc{\`a}-Ferretti system, \S5.4 \\
       & ? & other Matteini systems, \S5.4 \\
      \hline
    \end{tabular} \\
    \vspace*{5mm}
    Table 1: Lagrangian fibrations by Jacobians and Prym varieties
  \end{center}
\end{table}

In the first row, the Beauville-Mukai system degenerates to a compactification of the $\mathrm{GL}(n,\mathbb{C})$-Hitchin system. In the second and third rows we don't have as precise a relation between the compact and non-compact examples, though we have observed some analogies between the structures of the Debarre system and the $\mathrm{SL}(n,\mathbb{C})$-Hitchin system, and similarly for their dual fibrations. The compact counterparts of the $\mathrm{SO}(2n,\mathbb{C})$- and $\mathrm{SO}(2n+1,\mathbb{C})$-Hitchin systems are not clear at this time, and nor are the non-compact counterparts of the Markushevich-Tikhomirov system, the (various) Matteini systems, and the Arbarello-Sacc{\`a}-Ferretti system. Some of the gaps in the table in the non-compact column might be filled by various generalizations of Hitchin systems, e.g., those involving meromorphic Higgs bundles, irregular connections, parabolic and parahoric bundles, etc.

In this article we have focused on Lagrangian fibrations by Jacobians and Prym varieties, but these are not the only known examples of Lagrangian fibrations. 
\begin{itemize}
\item A very simple example comes from taking the Hilbert scheme of an elliptic K3 surface: this yields a Lagrangian fibration whose fibres are products of elliptic curves.
\item A Lagrangian fibration is {\em isotrivial\/} if all of its smooth fibres are isomorphic to a fixed abelian variety. Non-compact examples come from moduli spaces of Higgs bundles on elliptic curves, and were studied by Thaddeus~\cite{thaddeus01}. Compact examples also exist in all dimensions, and were studied by the author in~\cite{sawon14}.
\item There are also examples of Lagrangian fibrations whose fibres are intermediate Jacobians. Non-compact examples whose fibres are intermediate Jacobians of Calabi-Yau threefolds were constructed by Donagi and Markman~\cite{dm96i, dm96ii}. Compact examples arise by considering intermediate Jacobians of Fano threefolds, for example, families of Fano threefolds containing a fixed K3 surface (see Beauville~\cite{beauville02}, Iliev and Manivel~\cite{im07}, Hwang and Nagai~\cite{hn08}) or families of Fano threefolds contained in a fixed fourfold (see Markushevich and Tikhomirov~\cite{mt03}, Kuznetsov and Markushevich~\cite{km09}, Laza, Sacc{\`a}, and Voisin~\cite{lsv17}). Interestingly, in the last of these articles the authors use an identification of the intermediate Jacobians with certain Prym varieties.
\end{itemize}

Returning to Lagrangian fibrations by Jacobians and Prym varieties, Table 2 summarizes the known compact examples in four dimensions.

\begin{table}[h!]
\begin{center}
\begin{tabular}{|c|c|c|}
\hline
\textbf{Fibration} & \textbf{Polarization type} & \textbf{Singularities} \\
\hline\hline
Beauville-Mukai system & $(1,1)$ & smooth \\
\hline
Markushevich-Tikhomirov system & $(1,2)$ & $\mathbb{C}^4/\mathbb{Z}_2$ \\
\hline
Arbarello-Sacc{\`a}-Ferretti system & $(1,1)$ & $\mathbb{C}^4/\mathbb{Z}_2$ \\
\hline
Debarre system & $(1,3)$ & smooth \\
\hline
dual Debarre system & $(1,3)$ & $\mathbb{C}^4/\mathbb{Z}_3$ \\
\hline
\end{tabular} \\
    \vspace*{5mm}
    Table 2: Compact Lagrangian fibrations in dimension four
\end{center}
\end{table}

\begin{remark}
The polarization types of the fibres are $(1,d)$ where $d=1,2,3$. Could there be examples with larger values of $d$? Possibly, though in~\cite{sawon08ii} the author proved that $d\leq 1036$ under the natural hypothesis that the general singular fibres are semistable degenerations of abelian surfaces. Although this gives an upper bound, there is no reason to expect this bound to be sharp.
\end{remark}

\begin{remark}
There appears to be some correlation between the kinds of isolated singularities (of the form $\mathbb{C}^4/\mathbb{Z}_p$) that arise, and the polarization type $(1,d)$ of the fibres. A possible explanation for this is that the total space comes from a global quotient of the form $Y/\mathbb{Z}_p$, and this produces both the isolated singularities and fibres that are quotients of principally polarized abelian surfaces.
\end{remark}

\begin{flushleft}
Department of Mathematics\hfill sawon@email.unc.edu\\
University of North Carolina\hfill sawon.web.unc.edu\\
Chapel Hill NC 27599-3250\\
USA\\
\end{flushleft}

\end{document}